\documentclass[12pt,reqno,oneside]{amsart}
 \usepackage{verbatim,amsmath,amssymb,cite,xspace}
\usepackage{color,cite,graphicx}

\theoremstyle{plain}
\newtheorem{theorem}{Theorem}[section]
\newtheorem{claim}{Claim}[section]
\newtheorem{lemma}{Lemma}[section]

\theoremstyle{definition}

\theoremstyle{remark}

\usepackage{color}

\usepackage{amssymb}
\usepackage{amssymb, latexsym}
\usepackage{amsmath}
\usepackage{amscd}
\usepackage{enumerate}
\usepackage{amsfonts}
\usepackage{graphicx}
\usepackage{epsfig}

\newcommand{\HL}{\dot{H}^1\times L^2}
\newcommand{\OR}{\overrightarrow}
\newcommand{\E}{\mathcal{E}}

\newcommand{\dual}{2^{\ast}}
\newcommand{\Snorm}{L_t^{\frac{d+2}{d-2}}L_x^{2\frac{d+2}{d-2}}}

\newcommand{\spartial}{ {\slash\!\!\! \partial} }

\numberwithin{equation}{section} 

\begin{document}
\title[]{Soliton resolution along a sequence of times with dispersive error for type II singular solutions to focusing energy critical wave equation}
\author{Hao Jia}

\maketitle
\noindent
{\bf Abstract.}  
In this paper, we study the soliton resolution conjecture for Type II singular solutions $\OR{u}(t)$ to the focusing energy critical wave equation in $R^d\times [0,T_+)$, with $3\leq d\leq 5$. Suppose that $\OR{u}$ has a singularity at $(x,t)=(0,T_+)$, we show that along a sequence of times $t_n\uparrow T_+$ and in a neighborhood of $(0,T_+)$, $\OR{u}(t_n)$ can be decomposed as the sum of a regular part in the energy space, a finite combination of modulated solitons (translation, scaling and Lorentz transform of steady states) and a residue term which goes to zero in the Strichartz norm. In addition, the residue term is asymptotically radial and can concentrate energy only in a thin annulus near $|x|=T_+-t_n$. Our main tools include a Morawetz estimate, similar to the one used for wave maps, and the use of a virial identity to eliminate dispersive energy in $|x|<\lambda (T_+-t_n)$ for any $\lambda\in(0,1)$.

\begin{section}{Introduction}
Consider the energy critical wave equation:
\begin{equation}\label{eq:main}
\partial_{tt}u-\Delta u=|u|^{\dual-2}u,\,\,{\rm in}\,\, R^d\times (0,\infty)
\end{equation}
with initial data $(u_0,u_1)\in \HL(R^d)$, $3\leq d\leq 5$. We use the notation $\dual=\frac{2d}{d-2}$.
Equation (\ref{eq:main}) is invariant under space-time translations and the scaling for $\lambda>0$
\begin{eqnarray}
u(x,t)&\to&u_{\lambda}(x,t):=\lambda^{\frac{d}{2}-1}u(\lambda x,\lambda t),\\
(u_0(x),\,u_1(x))&\to&(u_{0\lambda}(x),\,u_{1\lambda}(x)):=(\lambda^{\frac{d}{2}-1}u_0(\lambda x),\,\lambda^{\frac{d}{2}}u_1(\lambda x)).
\end{eqnarray}
If $u$ is a solution to equation (\ref{eq:main}) with initial data $(u_0,u_1)$, then $u_{\lambda}$ is also a solution to equation (\ref{eq:main}) with initial data $(u_{0\lambda},\,u_{1\lambda})$.
A less obvious invariance for equation (\ref{eq:main}) is the Lorentz transform:
\begin{equation}
u(x,t)\,\to u_{\ell}(x,t):=u\left(x-\frac{x\cdot \ell}{|\ell |^2}\ell +\frac{\frac{x\cdot\ell}{|\ell |}\frac{\ell}{|\ell |}-\ell t}{\sqrt{1-|\ell |^2}},\,\frac{t-x\cdot\ell}{\sqrt{1-|\ell|^2}}\right),
\end{equation}
for each $\ell\in R^d$, with $|\ell |<1$. If $u$ is a solution to equation (\ref{eq:main}), then $u_{\ell}$ is also a solution where it is defined. 

\hspace{.6cm}The Cauchy problem for equation (\ref{eq:main}) has been well studied. For each $(u_0,u_1)\in\HL(R^d)$, there exists a unique solution $\OR{u}(t)=(u,\partial_tu)$ satisfying $\OR{u}\in C([0,T_+),\HL)$ and $u\in \Snorm(R^d\times [0,T))$ for any $T<T_+$. Moreover, the energy 
\begin{equation}
\E(\OR{u}(t)):=\int_{R^d}\,\left[\frac{|\partial_tu|^2}{2}+\frac{|\nabla u|^2}{2}-\frac{|u|^{\dual}}{\dual}\right](x,t)\,dx
\end{equation}
is preserved for $t\in [0,T_+)$. Here $T_+$ denotes the maximal time of existence for the solution $\OR{u}$. In addition, if $\|(u_0,u_1)\|_{\HL}$ is sufficiently small, then $T_+=\infty$ and the solution scatters, i.e., $u\in \Snorm(R^d\times [0,\infty))$. It is also well known that in general finite energy solutions to equation (\ref{eq:main}) may blow up in finite time, i.e., $T_+<\infty$. Indeed, using the finite speed of propagation for equation (\ref{eq:main}) to localize ODE type blow up solutions, one can easily construct finite time blow up solution $\OR{u}$ with $\|\OR{u}(t)\|_{\HL}\to\infty$, as $t\to T_+$. To rule out the ODE type behavior, one can focuse on the type II solutions. A solution $\OR{u}(t)$ is called Type II, if
\begin{equation*}
\sup_{t\in [0,T_+)}\,\|\OR{u}(t)\|_{\HL}<\infty.
\end{equation*}
Equation (\ref{eq:main}) admits infinitely many finite energy steady states $Q\in\dot{H}^1$, i.e., 
\begin{equation}
-\Delta Q=|Q|^{\dual-2}Q\,,\quad {\rm in}\,\,R^d.
\end{equation}
Among them, a distinguished role is played by the {\it ground state} $W$, which is the unique (up to scaling symmetry and sign change) radial $\dot{H}^1$ steady state. $W$ can be characterized as the minimizer of the Sobolev embedding $\dot{H}^1\subset L^{\dual}(R^d)$, see \cite{Talenti}.
We can use the symmetries of equation (\ref{eq:main}) and the steady states to obtain travelling wave solutions for $\ell\in R^d$ with $|\ell |<1$:
\begin{equation}
Q_{\ell}(x,t):=Q\left(x-\frac{x\cdot \ell}{|\ell |^2}\ell +\frac{\frac{x\cdot\ell}{|\ell |}\frac{\ell}{|\ell |}-\ell t}{\sqrt{1-|\ell |^2}}\right)
\end{equation}
$Q_{\ell}$ travels in the direction of $\ell$ with speed $|\ell |$. 

\hspace{.6cm}Type II solutions can have rich dynamics. Firstly, small solutions are Type II, exist globally and scatter. The travelling wave solutions $Q_{\ell}$ are global, but do not scatter. In addition, in \cite{HillRap,KSTwave}, Type II blow up solutions in the form of a rescaled ground state plus a small dispersive term  were constructed. More precisely the solution is given by 
\begin{equation*}
u(x,t)=\lambda(t)^{-\frac{d}{2}+1}W\left(\frac{x}{\lambda(t)}\right)+\epsilon(x,t),
\end{equation*}
where $\lambda(t)\to 0+$ as $t\to T_+$, and $\OR{\epsilon}(t)=(\epsilon,\,\partial_t\epsilon)$ is small in the energy space. It is expected that multi-soliton concentration is also possible for Type II blow up solutions, and it is an open problem to construct such a solution. Similar blow up solutions have been constructed for the energy critical equivariant wave maps and the radial energy critical Yang Mills equation, see \cite{KSTmills, KSTwavemap, RapRod, RodSter}.

\hspace{.6cm}The question of characterization of singularity formation or long time dynamics for type II solutions to equation (\ref{eq:main}) has been intensively studied. The central problem is to prove the ``soliton resolution conjecture". In the context of equation (\ref{eq:main}), the soliton resolution conjecture predicts that any type II solution should asymptotically decouple into a finite sum of modulated solitons, a regular part in the finite time blow up case or a free radiation in the global case, plus a residue term that vanishes asymptotically in the energy space as time approaches the maximal time of existence. The soliton resolution conjecture was proved in a remarkable work \cite{DKM} by Duyckaerts, Kenig and Merle, in the radial case for $d=3$. For other dimensions in the radial case, only soliton resolution along a sequence of times was available, see \cite{Kenig4dWave,Casey,JiaKenig}. In the nonradial case, the soliton resolution was proved in \cite{DKMsmallnonradial} for $d=3,\,5$, under an extra smallness condition. Without any smallness condition, perhaps the best result was \cite{DKMprofile}, where along a sequence of times, the smallest scale behavior of the solution was characterized as modulated solitons. 

\hspace{.6cm}In this paper, we consider the case of Type II blow up solution $\OR{u}(t)$, and prove the following result.
\begin{theorem}\label{th:main}
Let $\OR{u}\in C([0,T_+),\HL(R^d))$, with $u\in \Snorm(R^d\times [0,T))$ for any $T<T_+<\infty$, be a Type II blow up solution to equation (\ref{eq:main}). Here $T_+$ denotes the maximal existence time of $u$. Define the singular set 
\begin{equation}
\mathcal{S}:=\left\{x_{\ast}\in R^d:\,\|u\|_{\Snorm\left(B_{\epsilon}(x_{\ast})\times [T_+-\epsilon,\,T_+)\right)}=\infty,\quad{\rm for\,\,any\,\,}\epsilon>0\right\}.
\end{equation}
Then $\mathcal{S}$ is a set of finitely many points only. Let $x_{\ast}\in\mathcal{S}$ be a singular point. Then there exist integer $J_{\ast}\ge 1$, $r_{\ast}>0$, $\OR{v}\in \HL$, time sequence $t_n\uparrow T_+$,  scales $\lambda_n^j$ with $0<\lambda_n^j\ll T_+-t_n$, positions $c_n^j\in R^d$ satisfying $c_n^j\in B_{\beta (T_+-t_n)}(x_{\ast})$ for some $\beta\in(0,1)$ with $\ell_j=\lim\limits_{n\to\infty}\frac{c_n^j}{T_+-t_n}$ well defined, and travalling waves $Q_{\ell_j}$, for $1\leq j\leq J_{\ast}$, such that inside the ball $B_{r_{\ast}}$ we have
\begin{equation}\label{eq:maindecompositionhaha}
\OR{u}(t_n)=\OR{v}+\sum_{j=1}^{J_{\ast}}\,\left((\lambda_n^j)^{-\frac{d}{2}+1}\, Q_{\ell_j}\left(\frac{x-c_n^j}{\lambda_n^j},\,0\right),\,(\lambda_n^j)^{-\frac{d}{2}}\, \partial_tQ_{\ell_j}\left(\frac{x-c_n^j}{\lambda_n^j},\,0\right)\right)+(\epsilon_{0n},\,\epsilon_{1n}).
\end{equation}
Let $\OR{\epsilon}^L_n$ be the solution to linear wave equation with initial data $(\epsilon_{0n},\,\epsilon_{1n})\in\HL$, then the follwing vanishing condition holds
\begin{equation}\label{eq:refinedintro1}
\|\epsilon^L_n\|_{\Snorm(R^d\times R)}+\|( \epsilon_{0n},\,\epsilon_{1n})\|_{\HL( \{x\in R^d:\,|x|<\lambda (T_+-t_n)\,\,{\rm or}\,\,|x|>T_+-t_n\})}+\|\spartial\,\epsilon_{0n}\|_{L^2(R^d)}\to 0,
\end{equation}
as $n\to\infty$, for any $\lambda\in(0,1)$. In the above $\spartial$ denotes the tangential derivative.
In addition, the parameters $\lambda_n^j,\,c_n^j$ satisfy the pseudo-orthogonality condition
\begin{equation}
\frac{\lambda_n^j}{\lambda_n^{j'}}+\frac{\lambda_n^{j'}}{\lambda_n^j}+\frac{\left|c^j_n-c^{j'}_n\right|}{\lambda^j_n}\to\infty,
\end{equation}
as $n\to\infty$, for each $1\leq j\neq j'\leq J_{\ast}$.
\end{theorem}

\smallskip
\noindent
{\it Remark.} Theorem \ref{th:main} gives a partial proof of the soliton resolution conjecture for equation (\ref{eq:main}) along a sequence of times. Unfortunately our arguments do not rule out the dispersive energy completely, and it remains open to show the most natural vanishing condition $\|(\epsilon_{0n},\,\epsilon_{1n})\|_{\HL}\to 0$ as $n\to\infty$.\\

\hspace{.6cm}Our main new input in the proof of Theorem \ref{th:main} is the observation that the main part 
\begin{equation}\label{eq:mainpartofenergyflux}
\int_t^{T_+}\int_{|x|=T_+-s}\,\left[\frac{|\nabla u|^2}{2}+\frac{|\partial_tu|^2}{2}-\frac{x}{|x|}\cdot\nabla u\,\partial_tu\right](x,s)\,dxds
\end{equation}
of the energy flux
\begin{equation*}
\int_t^{T_+}\int_{|x|=T_+-s}\,\left[\frac{|\nabla u|^2}{2}+\frac{|\partial_tu|^2}{2}-\frac{x}{|x|}\cdot\nabla u\,\partial_tu-\frac{|u|^{\dual}}{\dual}\right](x,s)\,dxds
\end{equation*}
can be controlled, despite the fact that the energy flux may not be positive in general. The control on the main flux term (\ref{eq:mainpartofenergyflux}) allows us to use a Morawetz estimate, similar to the one used in the energy critical wave maps, see \cite{GrillakisEnergy,taoIII,Tataru4}. The Morawetz estimate is already sufficient to obtain a characterization of all profiles which are in some sense centered at time $t=0$, in the profile decomposition of $\OR{u}(t_n)$ along a well chosen time sequence $t_n$. As a result, we obtain a preliminary decomposition similar to the decomposition (\ref{eq:maindecompositionhaha}) but with the residue term vanishing in $L^{\dual}$ only.  To characterize other profiles that come from time infinity and obtain the refined vanishing condition (\ref{eq:refinedintro1}), we use a virial identity. The idea of using additional virial identities to eliminate dispersive energy was first introduced in \cite{JiaKenig}. In our case, an interesting new feature is that the virial identity can be useful only after we have obtained the preliminary decomposition, unlike the radial case in \cite{JiaKenig} where the virial identity gives good vanishing condition immediately. The difference is that in the radial case, we know that there is asymptotically no energy in the self similar region, see \cite{Kenig4dWave, JiaKenig}, while this is not true in general in the nonradial case. The use of an additional vanishing condition coming from the virial identity allows us to further characterize the residue term, and eventually prove (\ref{eq:refinedintro1}). We remark that the use of the virial identity appears to be crucial not only for ruling out dispersive energy in $B_{\lambda (T_+-t_n)}$ for $\lambda\in(0,1)$, but also for ruling out certain profiles coming from time infinity, and the Morawetz estimate, although sufficient for the purpose of characterzing profiles centered at $t=0$, seems not strong enough to ruling out certain other profiles, in our arguments.

\hspace{.6cm} Our paper is organized as follows. In section 2 we recall basic properties of profile decompostions for wave equations and properties of solutions to the linear wave equation, in section 3 we derive the Morawetz estimate, in section 4 we obtain some vanishing condition from the Morawetz estimate, in section 5 using the vanishing condition we obtain a preliminary decomposition, in section 6 we use the virial identity to derive further vanishing condition and in section 7 we conclude the proof.

\section{Preliminaries on Profile decompositions}
We briefly recall the profile decompositions, which was firstly introduced to wave equations by Bahouri and Gerard \cite{BaGe} in $R^3$, and was then extended to general dimensions by Bulut \cite{Bulut}. 

\hspace{.6cm}Let $(u_{0n},u_{1n})\in \HL$ satisfy $\sup\limits_{n}\|(u_{0n},u_{1n})\|_{\HL}\leq M<\infty$. Then passing to a subsequence there exist scales $\lambda^j_n>0$, positions $c^j_n\in R^d$, time translations $t^j_n\in R$, and finite energy solutions $U_j^L$ to the linear wave equation for each $j$, such that we have the following decomposition
\begin{eqnarray}
(u_{0n},\,u_{1n})&=&\sum_{j=1}^J\,\left(\left(\lambda^j_n\right)^{-\frac{d}{2}+1}U_j^L\left(\frac{x-c^j_n}{\lambda^j_n}, \,-\frac{t^j_n}{\lambda^j_n}\right),\,\left(\lambda^j_n\right)^{-\frac{d}{2}}\partial_tU_j^L\left(\frac{x-c^j_n}{\lambda^j_n}, \,-\frac{t^j_n}{\lambda^j_n}\right)\right)\nonumber\\
&&\nonumber\\
&&\quad\quad\quad+\left(w^J_{0n},\,w^J_{1n}\right),\label{eq:profiledecomposition}
\end{eqnarray}
where the parameters satisfy
\begin{equation}\label{eq:profileorthogonality1}
t^j_n\equiv 0,\,\,{\rm for\,\,all\,\,}n,\,\,{\rm or}\,\,\lim_{n\to\infty}\frac{t^j_n}{\lambda^j_n}\in\{\pm \infty\},
\end{equation}
and for $j\neq j'$
\begin{equation}\label{eq:profileorthogonality2}
\lim_{n\to\infty}\,\left(\frac{\lambda^j_n}{\lambda_n^{j'}}+\frac{\lambda^{j'}_n}{\lambda^j_n}+\frac{\left|c^j_n-c^{j'}_n\right|}{\lambda^j_n}+\frac{\left|t^j_n-t^{j'}_n\right|}{\lambda^j_n}\right)=\infty.
\end{equation}
In addition, let $\OR{w}^J_n$ be the solution to the linear wave equation with $\OR{w}^J_n(0)=\left(w^J_{0n},\,w^J_{1n}\right)$, then $w^J_n$ vanishes in the sense that
\begin{equation}
\lim_{J\to\infty}\,\limsup_{n\to\infty}\left(\|w^J_n\|_{\Snorm(R^d\times R)}+\|w^J_n\|_{L^{\infty}_tL^{\dual}_x(R^d\times R)}\right)=0.
\end{equation}
Moreover, if we write for $j\leq J$
\begin{equation}\label{eq:tildewj}
\OR{w}^J_n(x,t)=\left((\lambda^j_n)^{-\frac{d}{2}+1}\widetilde{w}^J_{jn}\left(\frac{x-c^j_n}{\lambda^j_n},\,\frac{t-t^j_n}{\lambda^j_n}\right),\,(\lambda^j_n)^{-\frac{d}{2}}\partial_t\widetilde{w}^J_{jn}\left(\frac{x-c^j_n}{\lambda^j_n},\,\frac{t-t^j_n}{\lambda^j_n}\right)\right),
\end{equation}
then $\OR{\widetilde{w}}^J_{jn}(t)=\left(\widetilde{w}^J_{jn},\,\partial_t\widetilde{w}^J_{jn}\right)\rightharpoonup 0$, as $n\to\infty$ for each $j\leq J$. \\

\hspace{.6cm}For later applications, we also need to recall some properties of free radiations (i.e, finite energy solutions to the linear wave equation) and of the profile decomposition. We begin with the following lemma which describes the concentration property of free radiations.
\begin{lemma}\label{lm:concentrationoffreeradiation}
Let $U$ be the solution to the $d+1$ dimensional linear wave equation in $R^d\times R$ with initial data $(U_0,\,U_1)\in\HL(R^d)$. Then
\begin{equation}\label{eq:concentrationforfreeradiation}
\lim_{\lambda\to\infty}\,\sup_{t\in R}\int_{\left||x|-|t|\right|>\lambda}\left[|\partial_tU|^2+|\nabla U|^2\right](x,t)\,dx=0.
\end{equation}
\end{lemma}

\smallskip
\noindent
{\it Proof.} For any $\epsilon>0$, there exists Schwartz class initial data $(u_0,\,u_1)$ such that
\begin{equation}
\|(u_0,\,u_1)-(U_0,\,U_1)\|_{\HL}<\epsilon.
\end{equation}
Let $\OR{u}$ be the solution to the $d+1$ dimensional linear wave equation with initial data $(u_0,\,u_1)$. Then by Theorem 1.1 in \cite{weightedestimate}, 
\begin{equation}\label{eq:boundforucon}
\sup_{(x,t)\in R^d\times R}(1+|x|+|t|)^{\frac{d-1}{2}}(1+||x|-|t||)^{\frac{d-1}{2}}|\nabla_{x,t}u(x,t)|<\infty.
\end{equation}
Hence for $u$, we can verify directly from (\ref{eq:boundforucon}) that 
\begin{equation}\label{eq:concentrationforu12}
\lim_{\lambda\to\infty}\,\sup_{t\in R}\int_{\left||x|-|t|\right|>\lambda}\left[|\partial_tu|^2+|\nabla u|^2\right](x,t)\,dx=0.
\end{equation}
Note that by energy conservation for linear wave equation, we have that
\begin{equation*}
\|\OR{u}(t)-\OR{U}(t)\|_{\HL}<\epsilon.
\end{equation*}
Thus,
\begin{equation}
\limsup_{\lambda\to\infty}\,\sup_{t\in R}\int_{\left||x|-|t|\right|>\lambda}\left[|\partial_tU|^2+|\nabla U|^2\right](x,t)\,dx\lesssim \epsilon^2.
\end{equation}
Since $\epsilon>0$ is arbitrary, (\ref{eq:concentrationforfreeradiation}) follows.\\

\hspace{.6cm}We shall also need the following lemma on the ``absolute continuity" of energy distribution for the free radiation on its concentration set $||x|-|t||=O(1)$, adapted to profiles. Denote
\begin{equation}
\OR{U}^L_{jn}(x,t)=\left(\left(\lambda^j_n\right)^{-\frac{d}{2}+1}U_j^L\left(\frac{x-c^j_n}{\lambda^j_n}, \,\frac{t-t^j_n}{\lambda^j_n}\right),\,\left(\lambda^j_n\right)^{-\frac{d}{2}}\partial_tU_j^L\left(\frac{x-c^j_n}{\lambda^j_n}, \,\frac{t-t^j_n}{\lambda^j_n}\right)\right).
\end{equation}
We have
\begin{lemma}\label{lm:absolutecontinuity}
Let $\OR{U}^L_{jn}(0)$ be a linear profile in $\HL$ with $\lim\limits_{n\to\infty}\frac{t^j_n}{\lambda^j_n}\in\{\pm\infty\}$. Suppose that a sequence of sets $E^j_n\subseteq R^d$ satisfy for all $M>0$ that
\begin{equation}
\left|\left\{x\in R^d:\,\left||x-c^j_n|-|t^j_n|\right|<M\lambda^j_n\right\}\cap\, E^j_n\right|=o\left(\left|t^j_n\right|^{d-1}\lambda^j_n\right),
\end{equation}
as $n\to\infty$. Then we have
\begin{equation}
\lim_{n\to\infty}\,\|\OR{U}^L_{jn}(0)\|_{\HL(E^j_n)}=0.
\end{equation}
\end{lemma}

\smallskip
\noindent
{\it Proof.} If $\OR{U}^L_{jn}$ has Schwartz class initial data, the lemma follows from the bound (\ref{eq:boundforucon}). The general case follows from approximation.\\

\hspace{.6cm}We shall also need the following orthogonality property for the profiles. 
We have
\begin{lemma}\label{eq:profileorthopre}
Suppose that $(u_{0n},\,u_{1n})$ is a bounded sequence in $\HL$ and has the profile decomposition (\ref{eq:profiledecomposition}). Then for each $J$
\begin{equation}\label{eq:almostequalno1}
\|(u_{0n},\,u_{1n})\|^2_{\HL}=\sum_{j=1}^J\|\OR{U}^L_j\|^2_{\HL}+\|w^J_{0n}\|_{\HL}^2+o_n(1),
\end{equation}
as $n\to\infty$, and for all $t$
\begin{eqnarray}
&&\int_{R^d}\,\left[(u_{1n}+\partial_1u_{0n})^2+|\nabla_{x'}u_{0n}|^2\right](x)\,dx\nonumber\\
&&\quad=\sum_{j=1}^J\int_{R^d}\,\left[\left(\partial_tU^L_j+\partial_1U^L_j\right)^2+|\nabla_{x'}U^L_j|^2\right](x,t)\,dx\nonumber\\
&&\quad\quad+\int_{R^d}\left[\left(w^J_n+\partial_1w^J_n\right)^2+|\nabla_{x'} w^J_n|^2\right](x,t)\,dx+o_n(1),\label{eq:almostequalno2}
\end{eqnarray}
as $n\to\infty$. In the above, $x=(x_1,x')$.
\end{lemma}

\smallskip
\noindent
{\it Proof.} The proof follows easily from the fact that for the linear wave $\OR{U}$, the quantities $\|\OR{U}\|^2_{\HL}$ and $\int\,(\partial_tU+\partial_1U)^2+|\nabla_{x'}U|^2\,dx$ are conserved, and the orthogonality condition for the profiles. Indeed, the first quantity is the energy, and the second quantity is a sum of the energy and the momentum in the $x_1$ direction. For the sake of completeness, we briefly outline the proof of (\ref{eq:almostequalno2}). It suffices to prove that for each $1\leq j\neq j'\leq J$\\
\begin{eqnarray}
\lim_{n\to\infty}\int_{R^d}\left[\left(\partial_t\,U^L_{jn}+\partial_1\,U^L_{jn}\right)\left(\partial_t\,U^L_{j'n}+\partial_1\,U^L_{j'n}\right)
+\nabla_{x'}\,U^L_{jn}\cdot\nabla_{x'}\,U^L_{j'n}\right](x,t)dx&=&0;\label{eq:orthogo11}\\
\lim_{n\to\infty}\int_{R^d}\left[\left(\partial_t\,U^L_{jn}+\partial_1\,U^L_{jn}\right)\left(\partial_t\,w^J_n+\partial_1\,w^J_n\right)
+\nabla_{x'}\,U^L_{jn}\cdot\nabla_{x'}\,w^J_n\right](x,t)dx&=&0.\label{eq:orthogo12}
\end{eqnarray}
It is easy to check that the left hand sides without taking limit are constant in time. Recalling the definition (\ref{eq:tildewj}) of $\OR{\widetilde{w}}^J_{jn}$ from the profile decompositions, we see by taking $t=t^j_n$ and rescaling that the left hand side of (\ref{eq:orthogo12}) is
\begin{equation*}
\int_{R^d}\,\left[\left(\partial_t\,U^L_j+\partial_1\,U^L_j\right)\left(\partial_t\,\widetilde{w}^J_{jn}+\partial_1\,\widetilde{w}^J_{jn}\right)
+\nabla_{x'}U^L_j\cdot\nabla_{x'}\widetilde{w}^J_{jn}\right](x,0)\,dx\to 0,
\end{equation*}
as $n\to\infty$, since $\widetilde{w}^J_{jn}(\cdot,t)\rightharpoonup 0$. The proof of (\ref{eq:orthogo11}) is similar. \\

\hspace{.6cm}The orthogonality conditions can be localized. More precisely we have
\begin{lemma}\label{lm:localizedorthogonality}
Suppose that $(u_{0n},\,u_{1n})$ is a bounded sequence in $\HL$ and has the profile decomposition (\ref{eq:profiledecomposition}). Fix $1\leq j\leq J$ and assume that
\begin{equation*}
\lim_{n\to\infty}\frac{t^j_n}{\lambda^j_n}\in\{\pm\infty\}.
\end{equation*}
Then for any $\beta>1$, 
\begin{equation}\label{eq:localizedenergyexpansion}
\liminf_{n\to\infty}\|(u_{0n},\,u_{1n})\|^2_{\HL\left(\left\{x:\,\frac{|t^j_n|}{\beta}<|x-c^j_n|<\beta |t^j_n|\right\}\right)}\ge \|\OR{U}^L_j\|_{\HL},
\end{equation}
and
\begin{eqnarray}
&&\liminf_{n\to\infty}\int_{\frac{|t^j_n|}{\beta}<|x-c^j_n|<\beta |t^j_n|}\,(u_{1n}+\partial_1u_{0n})^2+|\nabla_{x'}u_{0n}|^2\,dx\nonumber\\
&&\ge\int_{R^d}\,\left[\left(\partial_tU^L_j+\partial_1U^L_j\right)^2+|\nabla_{x'}U^L_j|^2\right](x,t)\,dx.\label{eq:localizedexpansion}
\end{eqnarray}
\end{lemma}

\smallskip
\noindent
{\it Proof.} Let us only prove (\ref{eq:localizedexpansion}), the proof of (\ref{eq:localizedenergyexpansion}) is similar. Let us assume $\lim\limits_{n\to\infty}\frac{t^j_n}{\lambda^j_n}=\infty$, the other case is identical. By Lemma \ref{lm:concentrationoffreeradiation} we get that
\begin{equation}\label{vanishing1111}
\lim_{n\to\infty}\int_{|x-c^j_n|>\beta t^j_n,\,\,{\rm or}\,\,|x-c^j_n|<\frac{t^j_n}{\beta}}\,|\nabla_{x,t}U^L_{jn}|^2(x,0)\,dx=0.
\end{equation}
Write 
\begin{equation}
(\widetilde{u}_{0n},\,\widetilde{u}_{1n})=(u_{0n},\,u_{1n})-\OR{U}^L_{jn}(\cdot,0).
\end{equation}
Then
\begin{eqnarray*}
&&\int_{\frac{t^j_n}{\beta}<|x-c^j_n|<\beta t^j_n}\,|u_{1n}+\partial_1u_{0n}|^2+|\nabla_{x'}u_{0n}|^2\,dx\\
&&=\int_{\frac{t^j_n}{\beta}<|x-c^j_n|<\beta t^j_n}\,\left[\left(\partial_tU^L_{jn}+\partial_1U^L_{jn}\right)^2+|\nabla_{x'}U^L_{jn}|^2\right](x,0)\,dx+\\
&&\quad\quad+\int_{\frac{t^j_n}{\beta}<|x-c^j_n|<\beta t^j_n}|\widetilde{u}_{1n}+\partial_1\widetilde{u}_{0n}|^2+|\nabla_{x'}\widetilde{u}_{0n}|^2\,dx+\\
&&\quad\quad\quad+\,2\int_{R^d}\,\left(\partial_tU^L_{jn}+\partial_1U^L_{jn}\right)\left(\widetilde{u}_{1n}+\partial_1\widetilde{u}_{0n}\right)+
\nabla_{x'}U^L_{jn}(x,0)\cdot\nabla_{x'}\widetilde{u}_{0n}\,dx+o_n(1)\\
&&\ge \int_{R^d}\left[\left(\partial_tU^L_{jn}+\partial_1U^L_{jn}\right)^2+|\nabla_{x'}U^L_{jn}|^2\right](x,0)\,dx+o_n(1).
\end{eqnarray*}
In the above, we have used the orthogonality between $\OR{U}^L_{jn}(x,0)$ and $(\widetilde{u}_{0n},\,\widetilde{u}_{1n})$. The lemma is proved.\\

For applications below, we shall also need the following version of the localized orthogonality with a small set removed.
\begin{lemma}\label{lm:removedcone}
Suppose that $(u_{0n},\,u_{1n})$ is a bounded sequence in $\HL$ and has the profile decomposition (\ref{eq:profiledecomposition}). Fix $1\leq j\leq J$ and assume that
\begin{equation*}
\lim_{n\to\infty}\frac{t^j_n}{\lambda^j_n}\in\{\pm\infty\}.
\end{equation*}
Assume that a sequence of sets $E^j_n\subseteq R^d$ satisfy for all $M>0$ that
\begin{equation}
\left|\left\{x\in R^d:\,\left||x-c^j_n|-|t^j_n|\right|<M\lambda^j_n\right\}\cap\, E^j_n\right|=o\left(\left|t^j_n\right|^{d-1}\lambda^j_n\right),
\end{equation}
as $n\to\infty$.
Then for any $\beta>1$, 
\begin{equation}\label{eq:coneenergyexpansion}
\liminf_{n\to\infty}\,\|(u_{0n},\,u_{1n})\|^2_{\HL\left(\left\{x:\,\frac{|t^j_n|}{\beta}<|x-c^j_n|<\beta |t^j_n|\right\}{\displaystyle \backslash}\,E^j_n\right)}\ge \|\OR{U}^L_j\|_{\HL}.
\end{equation}
\end{lemma}

\smallskip
\noindent
{\it Proof.} Assume that $\lim\limits_{n\to\infty}\frac{t^j_n}{\lambda^j_n}=\infty$, the other case is identical. Fix $\beta>1$, denote
\begin{equation}
F^j_n:=\left\{x:\,|x-c^j_n|>\beta\,t^j_n\,\,{\rm or}\,\,|x-c^j_n|<\frac{t^j_n}{\beta}\right\}\bigcup\,E^j_n.
\end{equation}
By Lemma \ref{lm:concentrationoffreeradiation} and Lemma \ref{lm:absolutecontinuity}, we see that
\begin{equation}
\left\|\OR{U}^L_{jn}(0)\right\|_{\HL\left(F^j_n\right)}\to 0,\,\,{\rm as}\,\,n\to\infty.
\end{equation}
Write 
\begin{equation}
(\widetilde{u}_{0n},\,\widetilde{u}_{1n})=(u_{0n},\,u_{1n})-\OR{U}^L_{jn}(\cdot,0).
\end{equation}
Then by Lemma \ref{lm:absolutecontinuity} and Lemma \ref{lm:localizedorthogonality}, we get that
\begin{eqnarray*}
&&\int_{\{x:\,\frac{t^j_n}{\beta}<|x-c^j_n|<\beta t^j_n\} \backslash\,E^j_n}\,|u_{1n}|^2+|\nabla u_{0n}|^2\,dx\\
&&=\int_{\{x:\,\frac{t^j_n}{\beta}<|x-c^j_n|<\beta t^j_n\}\backslash\,E^j_n}\,\left[\left(\partial_tU^L_{jn}\right)^2+|\nabla U^L_{jn}|^2\right](x,0)\,dx+\\
&&\quad\quad+\int_{\{x:\,\frac{t^j_n}{\beta}<|x-c^j_n|<\beta t^j_n\} \backslash\,E^j_n}|\widetilde{u}_{1n}|^2+|\nabla \widetilde{u}_{0n}|^2\,dx+\\
&&\quad\quad\quad+2\int_{R^d}\,\nabla U^L_{jn}(0)\,\nabla\widetilde{u}_{0n}+\partial_tU^L_{jn}(0)\,\partial_t\widetilde{u}_{1n}\,dx\\
&&\quad\quad\quad\quad\quad+O\left(\|\OR{U}^L_{jn}(0)\|_{\HL(F^j_n)}\right)\\
&&\ge \int_{R^d}\left[\left(\partial_tU^L_{jn}\right)^2+|\nabla_{x'}U^L_{jn}|^2\right](x,0)\,dx+o_n(1).
\end{eqnarray*}
In the above, we again used the pseudo-orthogonality of the profile $\OR{U}^L_{jn}$ and $(\widetilde{u}_{0n},\,\widetilde{u}_{1n})$. The lemma is proved.\\

\hspace{.6cm}To apply the linear profile decompositions to the nonlinear equation (\ref{eq:main}), we  need the following perturbation lemma.
\begin{lemma}\label{lm:longtimeperturbation}
Let $I$ be a time interval with $0\in I$. Let $\OR{U}\in C(I,\,\HL)$, with 
\begin{equation*}
\|U\|_{\Snorm(R^d\times I)}\leq M<\infty,
\end{equation*}
 be a solution to 
\begin{equation}\label{eq:equationwitherror}
\partial_{tt}U-\Delta U=|U|^{\dual-2}U+f,\,\,{\rm in}\,\,R^d\times I.
\end{equation}
Then we can find $\epsilon_{\ast}=\epsilon_{\ast}(M,d)>0$ sufficiently small, such that if
\begin{equation}\label{eq:errorissmall}
\|f\|_{L^1_tL^2_x(R^d\times I)}+\|(u_0,u_1)-\OR{U}(0)\|_{\HL}\leq \epsilon\leq \epsilon_{\ast},
\end{equation}
then the unique solution $\OR{u}$ to equation (\ref{eq:main}) with initial data $(u_0,u_1)\in\HL$ exists in $R^d\times I$. Moreover, $\OR{u}$ verifies the estimate
\begin{equation}
\sup_{t\in I}\,\|\OR{u}(t)-\OR{U}(t)\|_{\HL}+\|u-U\|_{\Snorm(R^d\times I)}\leq C(M,d)\epsilon.
\end{equation} 
\end{lemma}

\smallskip
\noindent
{\it Remark.} This perturbation lemma is often termed ``long time perturbation" lemma, due to the fact there is no restriction on the size of the time interval $I$. Such results are well known, see for example Theorem 2.14 in \cite{DKMacta} and Lemma 2.1 in \cite{JiaLiuXu}. We give a proof for the sake of completeness. \\

\smallskip
\noindent
{\it Proof.} By splitting the interval $I$ and the time reversibility of equation (\ref{eq:main}), we can assume without loss of generality that $I=[0,T]$. Fix $\beta=\beta(d)>0$ small, whose value is to be determined below. We divide the interval $I$ into $\bigcup_{j=0}^{N-1}\,[t_j,\,t_{j+1}]$, such that $t_0=0$, $t_N=T$, and that 
\begin{equation}
\|U\|_{\Snorm(R^d\times[t_j,\,t_{j+1}])}\leq \beta.
\end{equation}
Clearly we shall need $N=N(M,d)$ intervals. Fix a $K(d)>1$ sufficiently large whose value is to be determined below, we shall prove by induction that 
\begin{equation}\label{eq:allbounds}
\sup_{t\in[t_j,\,t_{j+1}]}\,\|\OR{u}(t)-\OR{U}(t)\|_{\HL}+\|u-U\|_{\Snorm(R^d\times[t_j,\,t_{j+1}])}\leq K^{j+1}\epsilon,
\end{equation}
for each $0\leq j\leq N-1$. Let $\omega=u-U$, then $\omega$ satisfies the equation
\begin{equation}\label{eq:perturbedequation7}
\partial_{tt}\omega-\Delta\omega=|U+\omega|^{\dual-2}(U+\omega)-|U|^{\dual-2}U-f,
\end{equation}
with initial data $\OR{\omega}(0)=(u_0,u_1)-\OR{U}(0)$. Let us note that the nonlinear term
\begin{equation}
\mathcal{N}(\omega):=|U+\omega|^{\dual-2}(U+\omega)-|U|^{\dual-2}U
\end{equation}
verifies 
\begin{equation}\label{eq:inequalitynonlinearterm}
\left|\mathcal{N}(\omega)\right|\lesssim_d\,|U|^{\dual-2}|\omega|+|\omega|^{\dual-1}.
\end{equation}
For $j\ge 0$, assume that 
\begin{equation}\label{eq:inductionstep}
\|\OR{\omega}(t_j)\|_{\HL}\leq K^j\epsilon,
\end{equation}
we shall establish the bound (\ref{eq:allbounds}) for the interval $[t_j,\,t_{j+1}]$. By local existence theory, for some $\tau>0$, we have 
\begin{equation}\label{eq:continuitybound}
\sup_{t\in[t_j,\,t_j+\tau]}\,\|\OR{\omega}(t)\|_{\HL}+\|\omega\|_{\Snorm(R^d\times[t_j,\,t_j+\tau])}\leq K^{j+1}\epsilon.
\end{equation}
We shall show that as long as the bound (\ref{eq:continuitybound}) holds and $K$ is chosen large, we actually have the improved bound
\begin{equation}
\sup_{t\in[t_j,\,t_j+\tau]}\,\|\OR{\omega}(t)\|_{\HL}+\|\omega\|_{\Snorm(R^d\times[t_j,\,t_j+\tau])}\leq \frac{1}{2}K^{j+1}\epsilon.
\end{equation}
Then by standard continuity arguments, we will be done. By the inequality (\ref{eq:inequalitynonlinearterm}) and the bound (\ref{eq:inductionstep}), we can estimate the nonlinear term as
\begin{equation*}
\|\mathcal{N}(\omega)\|_{L^1_tL^2_x(R^d\times[t_j,\,t_j+\tau))}\lesssim_d \beta^{\dual-2}K^{j+1}\epsilon+K^{(\dual-1)(j+1)}\epsilon^{\dual-1}.
\end{equation*}
Using the bound (\ref{eq:inductionstep}) and Strichartz estimates, we then get that
\begin{eqnarray}
&&\sup_{t\in[t_j,\,t_j+\tau]}\,\|\OR{\omega}(t)\|_{\HL}+\|\omega\|_{\Snorm(R^d\times[t_j,\,t_j+\tau])}\nonumber\\
&&\quad\leq\, C(d) K^j\epsilon+ C(d)\beta^{\dual-2}K^{j+1}\epsilon+C(d)K^{(\dual-1)(j+1)}\epsilon^{\dual-1}.\label{eq:boundfromstrichartz}
\end{eqnarray}
Now we choose $\beta=\beta(d)>0$ such that 
\begin{equation*}
C(d)\beta^{\dual-2}\leq\frac{1}{8}.
\end{equation*}
Then we fix $K=K(d)>0$ such that
\begin{equation*}
C(d)\leq \frac{K}{8}.
\end{equation*}
Then we fix $\epsilon_{\ast}(M,d)>0$ such that
\begin{equation*}
C(d)K^{(\dual-1)N}\epsilon_{\ast}<\frac{1}{8}.
\end{equation*}
With these choices of parameters, we can then conclude from (\ref{eq:boundfromstrichartz}) that
\begin{equation*}
\sup_{t\in[t_j,\,t_j+\tau]}\,\|\OR{u}(t)\|_{\HL}+\|u\|_{\Snorm(R^d\times[t_j,\,t_j+\tau])}\leq \frac{1}{2}K^{j+1}\epsilon.
\end{equation*}
The lemma is proved.\\

The perturbation lemma \ref{lm:longtimeperturbation} has important applications to the nonlinear profile decompositions. Assume that the sequence of initial data $(u_{0n},\,u_{1n})$ is uniformly bounded in the energy space with respect to $n$ and that $(u_{0n},\,u_{1n})$ has the profile decomposition (\ref{eq:profiledecomposition}). For each $j$, introduce the nonlinear profile $U_j$ as follows. 
\begin{itemize}
\item if $\lim\limits_{n\to\infty}\frac{t^j_n}{\lambda^j_n}=L\in\{\pm\infty\}$, then define $\OR{U}_j$ as the unique solution to equation (\ref{eq:main}) in a neighborhood of $-L$, with 
\begin{equation*}
\lim_{t\to -L}\|\OR{U}_j(t)-\OR{U}_j^L(t)\|_{\HL}=0.
\end{equation*}
The existence of $\OR{U}_j$ follows from standard perturbative arguments.
\item if $t^j_n\equiv 0$ for all $n$, then define $\OR{U}_j$ as the solution to equation (\ref{eq:main}) with initial data $\OR{U}_j(0)=\OR{U}_j^L(0)$.\\
\end{itemize}

We have the following nonlinear approximation lemma for the linear profile decomposition, see Bahouri and Gerard \cite{BaGe} for details.
\begin{lemma}\label{lm:nonlinearprofiledecomposition}
Let $(u_{0n},\,u_{1n})$ be a sequence of initial data that are uniformly bounded in the energy space. Assume that $(u_{0n},\,u_{1n})$ has the profile decomposition (\ref{eq:profiledecomposition}). Let $\OR{U}_j$ be the nonlinear profile associated with $\OR{U}_j^L$, $\lambda^j_n$, $t^j_n$. Denote $T_+(U_j)$ as the maximal time of existence for $U_j$. We adopt the convention that $T_+(U_j)=\infty$ if $\lim\limits_{n\to\infty}\frac{-t^j_n}{\lambda^j_n}=\infty$, and that the maximal interval of existence of $\OR{U}_j$ is $(-\infty,\,T_+(U_j))$ if $\lim\limits_{n\to\infty}\frac{-t^j_n}{\lambda^j_n}=-\infty$. Let $\theta_n\in (0,\infty)$. Assume that $\forall\, j$ and large $n$
\begin{equation}
\frac{\theta_n-t^j_n}{\lambda^j_n}<T_+(U_j),\,\,{\rm for}\,\,\forall \,n,\,\,{\rm and}\,\,\limsup_{n\to\infty}\left\|U^j\right\|_{\Snorm\left(R^d\times \left(-\frac{t^j_n}{\lambda^j_n},\,\frac{\theta_n-t^j_n}{\lambda^j_n}\right)\right)}<\infty.
\end{equation}
Let $\OR{u}_n$ be the solution to equation (\ref{eq:main}) with initial data $(u_{0n},\,u_{1n})$.
Then for sufficiently large $n$, $\OR{u}_n$ has the following decomposition
\begin{equation}
\OR{u}_n=\sum_{j=1}^J\,\left(\left(\lambda^j_n\right)^{-\frac{d}{2}+1}U_j\left(\frac{x-c^j_n}{\lambda^j_n}, \,\frac{t-t^j_n}{\lambda^j_n}\right),\,\left(\lambda^j_n\right)^{-\frac{d}{2}}\partial_tU_j\left(\frac{x-c^j_n}{\lambda^j_n}, \,\frac{t-t^j_n}{\lambda^j_n}\right)\right)+\OR{w}^J_n+\OR{r}^J_n,
\end{equation}
in $R^d\times [0,\,\theta_n)$, where $\OR{r}_n$ vanishes in the sense that
\begin{equation}
\lim_{J\to\infty}\,\limsup_{n\to\infty}\left(\sup_{t\in[0,\theta_n)}\|\OR{r}^J_n(t)\|_{\HL}+\|r^J_n\|_{\Snorm(R^d\times [0,\,\theta_n))}\right)=0.
\end{equation}
\end{lemma}

\medskip

The principle of finite speed of propagation plays an essential role in the study of wave equations. Below we shall use the following version of this principle. Let us set $a\wedge b:=\min\{a,\,b\}$ for any $a,\,b\in R$.
\begin{lemma}
Let $\OR{u},\,\OR{v}\in C([-T,T],\,\HL)$, with $u,\,v\in \Snorm(R^d\times[-T,T])$, be two solutions to equation (\ref{eq:main}), with initial data $(u_0,\,u_1)\in\HL$ and $(v_0,\,v_1)\in\HL$ respectively. Assume that for some $R>0$,
\begin{equation*}
(u_0(x),\,u_1(x))=(v_0(x),\,v_1(x))\,\,\,{\rm in}\,\,B_R.
\end{equation*}
Then 
\begin{equation}
u(x,t)=v(x,t)\,\,\,{\rm for}\,\,(x,t)\in\{(x,t):\,|x|<R-|t|,\,\,|t|<R\wedge T\}.
\end{equation}
\end{lemma}

The perturbation lemma \ref{lm:longtimeperturbation} combined with finite speed of propagation implies the following local-in-space approximation result.
\begin{lemma}\label{lm:localinspaceapproximation}
Let $\OR{u}_n$ and $\OR{v}_n$ be two sequences of solutions defined in $R^d\times [-T,T]$ with initial data $(u_{0n},\,u_{1n})\in\HL$ and $(v_{0n},\,v_{1n})\in\HL$ respectively. Assume that 
\begin{equation*}
\|u_n\|_{\Snorm(R^d\times [-T,T])}\leq M<\infty, 
\end{equation*}
and that 
\begin{equation*}
\|(u_{0n},\,u_{1n})-(v_{0n},\,v_{1n})\|_{\HL(B_{R+\epsilon})}+\|u_{0n}-v_{0n}\|_{L^{\dual}(B_{R+\epsilon})}\to 0,
\end{equation*}
as $n\to\infty$, for some $R>0$ and $\epsilon>0$ . Then we have for $\tau=R\wedge T$ that
\begin{equation}
\sup_{t\in(-\tau,\,\tau)}\,\|\OR{u}_n(t)-\OR{v}_n(t)\|_{\HL(|x|<R-|t|)}+\|u_n-v_n\|_{\Snorm(\{(x,t):\,|x|<R-|t|,\,|t|<\tau\})}\to 0,
\end{equation}
as $n\to\infty$.
\end{lemma}  

\smallskip
\noindent
{\it Proof.} We define the following modified initial data $(\tilde{v}_{0n},\,\tilde{v}_{1n})$, with
\begin{equation}
(\tilde{v}_{0n},\,\tilde{v}_{1n}):=\left\{\begin{array}{ll}
                              (u_{0n},\,u_{1n})&\,\,{\rm if}\,\,|x|>R+\epsilon,\\
                   &\\
                               \frac{|x|-R}{\epsilon}(u_{0n},\,u_{1n})+\frac{R+\epsilon-|x|}{\epsilon}(v_{0n},\,v_{1n}),&\,\,{\rm if}\,\,R<|x|<R+\epsilon,\\
&\\
                      (v_{0n},\,v_{1n}) &\,\,{\rm if}\,\,|x|<R.
                              \end{array}\right.
\end{equation}
Then it is straightforward to verify that
\begin{equation*}
\|(\tilde{v}_{0n},\,\tilde{v}_{1n})-(u_{0n},\,u_{1n})\|_{\HL(R^d)}\to 0,\,\,{\rm as}\,\,n\to\infty,
\end{equation*}
and that
\begin{equation*}
(v_{0n},\,v_{1n})\equiv (\tilde{v}_{0n},\,\tilde{v}_{1n}),\,\,\,{\rm in}\,\,B_R.
\end{equation*}
Let $\OR{\tilde{v}}_n$ be the solution to equation (\ref{eq:main}) with initial data $\OR{\tilde{v}}_n(0)=(\tilde{v}_{0n},\,\tilde{v}_{1n})$. Then by the principle of finite speed of propagation, we see that 
\begin{equation}\label{eq:localequality}
v_n\equiv \tilde{v}_n, \,\,\,\,{\rm in}\,\, \{(x,t):\,|x|<R-|t|,\,|t|<\tau\}.
\end{equation}
By the perturbation Lemma \ref{lm:longtimeperturbation}, we get that
\begin{equation}\label{eq:differencegoestozero}
\lim_{n\to\infty}\,\left(\sup_{t\in[-T,\,T]}\|\OR{u}_n(t)-\OR{\tilde{v}}_n(t)\|_{\HL(R^d)}+\|u_n-\tilde{v}_n\|_{\Snorm(R^d\times [-T,T])}\right)=0.
\end{equation}
Combining (\ref{eq:differencegoestozero}) with (\ref{eq:localequality}) finishes the proof of the lemma.\\

\end{section}

\begin{section}{Control of energy flux and Morawetz identity}
In this section we begin the study of the type II singular solutions to equation (\ref{eq:main}). Assume without loss of generality (by a rescaling argument) that $\OR{u}\in C((0,1],\HL)$, with $u\in \Snorm(R^d\times (\epsilon,1])$ for all $\epsilon>0$, is a Type II solution to equation (\ref{eq:main}) in $R^d\times (0,1]$ with Cauchy data prescribed at time $t=1$, and that $t=0$ is the blow up time for $u$, i.e., $\|u\|_{\Snorm(R^d\times (0,1])}=\infty$. This choice of convention saves us some trouble with notations below. Define the set of singular points
\begin{equation}
\mathcal{S}:=\left\{x_{\ast}\in R^d:\,\|u\|_{\Snorm(B_{\epsilon}(x_{\ast})\times (0,\epsilon])}=\infty,\,\,{\rm for\,\,all\,\,}\epsilon>0\right\}.
\end{equation}
By small data global existence theory for equation (\ref{eq:main}) and finite speed of propagation, we see that if $x_{\ast}\in\mathcal{S}$, then there exists a $\epsilon_0(d)>0$, with
\begin{equation}\label{eq:nontrivialcon}
\forall \,\epsilon>0,\,\,\liminf_{t\to 0+}\,\int_{|x-x_{\ast}|<\epsilon}\,\left[\frac{|\nabla u|^2}{2}+\frac{|\partial_tu|^2}{2}+\frac{|u|^{\dual}}{\dual}\right](x,t)\,dx>\epsilon_0.
\end{equation}
By the assumption that $\OR{u}$ is Type II, i.e., $\sup\limits_{t\in(0,1]}\,\|\OR{u}(t)\|_{\HL}<\infty$, Sobolev embedding $\dot{H}^1(R^d)\subset L^{\dual}(R^d)$ and the concentration of energy near singular points (\ref{eq:nontrivialcon}), we conclude that $\mathcal{S}$ is a non-empty set of finitely many points. Below we shall assume that $0\in \mathcal{S}$ and the main focus is to study the asymptotics of $u$ near this singular point. The finite speed of propagation, together with local Cauchy theory for equation (\ref{eq:main}), implies that 
$$\|u\|_{\Snorm\left((R^d\backslash B_R)\times (0,1]\right)}<\infty$$
for some sufficiently large $R>1$. Hence via a compactness argument, we see that
\begin{equation}
\|u\|_{\Snorm\left((R^d\backslash \bigcup_{j=1}^{K}B_{\epsilon}(x_j))\times (0,1]\right)}\leq C(\epsilon)<\infty,
\end{equation} 
where $\mathcal{S}=\{x_1,\dots,x_{K}\}$. In addition, $\OR{u}\in C([0,1],\,\HL(R^d\backslash \cup_{j=1}^{K}B_{\epsilon}(x_j)))$ for any $\epsilon>0$. Thus 
\begin{equation}
(f,g):=w-\lim_{t\to 0+}\,\OR{u}(t)
\end{equation}
is well defined as a function in $\HL(R^d)$. 

\hspace{.6cm}Let $\OR{v}\in C([0,\delta],\,\HL)$ with $v\in \Snorm(R^d\times [0,\delta])$ be the local in time solution to equation (\ref{eq:main}) with initial data $\OR{v}(0)=(f,g)$. Using finite speed of propagation again,  we have that
\begin{equation}\label{eq:regularsingular}
u(x,t)=v(x,t)\,\,{\rm for}\,\,(x,t)\in R^d\times (0,\delta]\backslash \bigcup_{j=1}^{K}\{(x,t):\,|x-x_j|\leq t, \,t\in (0,1]\}.
\end{equation}
$v$ is called the regular part of $u$. (\ref{eq:regularsingular}) shows that there is a sharp division of the singular and regular regions for type II solutions.

\hspace{.6cm}We now make the following simple observation. Since $v\in \Snorm(R^d\times [0,\delta])$ and $\OR{v}\in L^{\infty}([0,\delta],\,\HL)$, we see that 
\begin{equation}
\left|\nabla_{x,t}\,|v|^{\dual}\right|\lesssim_d |\nabla_{x,t} v||v|^{\dual-1}\in L^1(R^d\times[0,\delta]).
\end{equation}
This regularity property implies that $v$ is well defined on any sufficiently regular hypersurfaces as $L^{\dual}$ function. Similar arguments can be made for $u$ in the region $R^d\times (\epsilon,\,1]$ for any $\epsilon>0$.
In particular we conclude that $u\equiv v\in L^{\dual}(\{(x,t):\,|x|=t,\,t\in(0,\delta)\})$ is well defined. 

This observation, although simple, opens up the possibility of using more monotonicity quantities than previously known for equation (\ref{eq:main}). One obstacle in the study of equation (\ref{eq:main}) is that the energy flux is in general not positive. The main issue is the term $|u|^{\dual}$ which comes with a minus sign in the energy flux. However, since $|u|^{\dual}$ is integrable on the boundary of the lightcone, this negative term can controlled. In the end we can then control {\it the main flux term}
\begin{equation}\label{eq:mainflux}
\int_{t_1}^{t_2}\int_{|x|=t}\,\frac{|\nabla_{x,t}u|^2}{2}+\partial_tu\,\frac{x}{|x|}\cdot\nabla u\,d\sigma dt.
\end{equation}
The control on the flux term (\ref{eq:mainflux}) is very useful. In particular we can now use a variant of Morawetz identity from the study of energy critical wave maps \cite{GrillakisEnergy,taoIII,Tataru4}. We turn to the details.
The energy flux identity in our setting is
\begin{eqnarray}
&&\frac{1}{\sqrt{2}}\int_{t_1}^{t_2}\int_{|x|=t}\,\frac{|\nabla_{x,t}u|^2}{2}+\partial_tu\,\frac{x}{|x|}\cdot\nabla u-\frac{|u|^{\dual}}{\dual}\,d\sigma dt\nonumber\\
&&=\int_{|x|\leq t_2}\left[\frac{|\nabla_{x,t}u|^2}{2}-\frac{|u|^{\dual}}{\dual}\right](x,t_2)\,dx-\int_{|x|\leq t_1}\left[\frac{|\nabla_{x,t}u|^2}{2}-\frac{|u|^{\dual}}{\dual}\right](x,t_1)\,dx,\label{eq:energyfluxidentity}
\end{eqnarray}
for $0<t_1<t_2<\delta$. Since a priori, $|\nabla_{x,t}u|$ is not classically defined on the boundary of the lightcone, the term $\int_{t_1}^{t_2}\int_{|x|=t}\,\frac{|\nabla_{x,t}u|^2}{2}+\partial_tu\,\frac{x}{|x|}\cdot\nabla u\,d\sigma dt$ requires an appropriate interpretation. Note however that the term $\int_{t_1}^{t_2}\int_{|x|=t}\,\frac{|u|^{\dual}}{\dual}\,d\sigma dt$ is well defined by discussions above. Hence by approximation by smooth solutions, 
\begin{eqnarray*}
\frac{1}{\sqrt{2}}f(t_1,t_2)&:=&\int_{|x|\leq t_2}\left[\frac{|\nabla_{x,t}u|^2}{2}-\frac{|u|^{\dual}}{\dual}\right](x,t_2)\,dx-\int_{|x|\leq t_1}\left[\frac{|\nabla_{x,t}u|^2}{2}-\frac{|u|^{\dual}}{\dual}\right](x,t_1)\,dx\\
&&\quad\quad\quad\quad\quad\quad+\,\frac{1}{\sqrt{2}}\int_{t_1}^{t_2}\int_{|x|=t}\frac{|u|^{\dual}}{\dual}\,d\sigma\,dt
\end{eqnarray*}
is well defined. This function can be interpreted as the main flux term.
We claim that
\begin{claim}\label{cl:realanalysis}
\begin{equation}\label{eq:limitofconcentration}
\lim_{t\to 0+}\,\int_{|x|<t}\,\frac{|\nabla_{x,t}u|^2}{2}-\frac{|u|^{\dual}}{\dual}(x,t)\,dx\,\,\,{\rm exists},
\end{equation}
and 
\begin{equation}\label{eq:vanishingflux}
\lim_{\tau\to 0+}\int_{\tau}^t\int_{|x|=s}\,\frac{|\spartial u|^2}{2}+\frac{\left|\partial_tu+\frac{x}{|x|}\cdot\nabla u\right|^2}{2}\,d\sigma ds\leq C(u),
\end{equation}
for some $C(u)$ independent of $t\in[0,\delta]$,
where $|\spartial u|^2$ denotes the tangential derivative $|\nabla u|^2-\left|\frac{x}{|x|}\cdot\nabla u\right|^2$. 
\end{claim}
The proof of Claim \ref{cl:realanalysis} is straightforward. The claimed bound (\ref{eq:vanishingflux}) is a direct consequence of the assumption that $\OR{u}$ is type II, $u\in L^{\dual}(\{|x|=t,\,t\in(0,\delta)\})$ and the energy flux identity.  (\ref{eq:limitofconcentration}) follows easily from the energy flux identity.\\

Now we can state and prove the crucial Morawetz estimates.
\begin{lemma}\label{lm:morawetz}
Let $u$ be as above. Then for any $0<10\,t_1<t_2<\delta$, we have
\begin{equation}
\int_{t_1}^{t_2}\int_{|x|<t}\left(\partial_tu+\frac{x}{t}\cdot\nabla u+\left(\frac{d}{2}-1\right)\frac{u}{t}\right)^2\,dx\,\frac{dt}{t}\leq C\left(\log{\frac{t_2}{t_1}}\right)^{\frac{d}{d+1}},
\end{equation}
for some $C$ independent of $t_1,\,t_2$.
\end{lemma}

\smallskip
\noindent
{\it Remark.} Estimates of this type have played an important role in the study of energy critical wave maps, see Grillakis \cite{GrillakisEnergy}, Tao \cite{taoIII}, and Sterbenz-Tataru \cite{Tataru4}. Our calculations below are essentially the same as those in \cite{taoIII}. We will however provide a detailed proof, due to an extra term $u$ coming from the scaling in our case, and the importance of the estimate in our analysis below.\\

\smallskip
\noindent
{\it Proof.} In the calculations below, we can initially assume that $u$ is smooth by standard approximation. Note that it is in general not possible to achieve an approximation up to the singular time $t=0$. However since we are interested in the estimates for the fixed $t_1,\,t_2$, such approximation is harmless, as long as the bound only depend on the bounds on $u$, such as the bound on the energy flux of $u$, but not on any additional smoothness assumptions. 

\hspace{.6cm}For $0<10\,t_1<t_2<\delta$, fix $\epsilon>0$ to be determined below. Set
\begin{equation}
\rho=\rho_{\epsilon}=\left((1+\epsilon^2)\,t^2-|x|^2\right)^{-\frac{1}{2}},\quad{\rm for}\,\,|x|<t.
\end{equation}
Set $x^{\alpha}=t$ if $\alpha=0$, and $x^{j}=x_j$ if $1\leq j\leq d$. Let $\partial^{\alpha}=-\partial_{\alpha}=-\partial_t$ if $\alpha=0$, and $\partial^{j}=\partial_j$ for $1\leq j\leq d$. We note that
\begin{equation*}
x^{\alpha}\partial_{\alpha}\rho=-\rho^3\,(1+\epsilon^2)\,t^2+x^j\rho^3\,x_j=-\rho.
\end{equation*}
We can rewrite equation (\ref{eq:main}) as
\begin{equation}\label{eq:variantmain}
\partial^{\alpha}\partial_{\alpha}u+|u|^{\dual-2}u=0.
\end{equation}
Multiplying equation (\ref{eq:variantmain}) with $\rho\,(x^{\alpha}\partial_{\alpha}u+\left(\frac{d}{2}-1)u\right)$, and integrating in the region $\{(x,t):\,|x|<t,\,t\in(t_1,t_2)\}$, we get that
\begin{eqnarray*}
0&=&\int_{t_1}^{t_2}\int_{|x|<t}\,\left(\partial^{\alpha}\partial_{\alpha}u+|u|^{\dual-2}u\right)\,\rho \,\left(x^{\beta}\partial_{\beta}u+\left(\frac{d}{2}-1\right)u\right)\,dxdt\\
&=&\int_{t_1}^{t_2}\int_{|x|<t}\,\rho\, x^{\beta}(\partial^{\alpha}\partial_{\alpha}u)\,\partial_{\beta}u+\rho\, x^{\beta}\partial_{\beta}\frac{|u|^{\dual}}{\dual}+\left(\frac{d}{2}-1\right)\rho\,(\partial^{\alpha}\partial_{\alpha}u)\,u+\left(\frac{d}{2}-1\right)\rho\, |u|^{\dual}\,dxdt\\
&=&\int_{t_1}^{t_2}\int_{|x|<t}\,\rho\, x^{\beta}\partial^{\alpha}\left(\partial_{\alpha}u\partial_{\beta}u\right)-\rho\, x^{\beta}\partial_{\beta}\frac{\partial_{\alpha}u\partial^{\alpha }u}{2}+\rho\, x^{\beta}\partial_{\beta}\frac{|u|^{\dual}}{\dual}+\\
&&\quad\quad+\left(\frac{d}{2}-1\right)\rho\,\partial^{\alpha}\left(\partial_{\alpha}u\,u\right)-\left(\frac{d}{2}-1\right)\rho\,\partial_{\alpha}u\partial^{\alpha}u+\left(\frac{d}{2}-1\right)\rho\,|u|^{\dual}\,dxdt,\\
&=&{\rm boundary\,\,terms\,\,B}\,+{\rm interior \,\,terms\,\,I}.
\end{eqnarray*}
In the above we use integration by parts to obtain the boundary terms B and interior terms I. We shall treat B and I separately. Let us firslty record the following bounds for $\rho$ in the region $\{(x,t):\,|x|<t,\,t\in (0,\delta)\}$:
\begin{equation}\label{eq:ineqrho}
|\rho(x,t)|\leq \epsilon^{-1}t^{-1}.
\end{equation}
Note also that the bound on the energy flux (\ref{eq:vanishingflux}) implies that
\begin{equation}\label{eq:usebound1}
\int_{t_1}^{t_2}\int_{|x|=t}\,\left|\frac{x^{\alpha}}{t}\partial_{\alpha} u\right|^2\,d\sigma dt\leq C(u)<\infty
\end{equation}
for some constant $C(u)$ independent of $t_1,\,t_2$. Note that the unit outer normal on the lateral boundary of the lightcone $|x|=t$ is $\overrightarrow{n}=(n_{\alpha})=\frac{1}{\sqrt{2}}\left(-1,\frac{x}{|x|}\right)$, and that 
\begin{equation}
x^{\beta}\,n_{\beta}\equiv 0,\,\,{\rm on}\,\,|x|=t.
\end{equation}
The boundary term $B$ can be written as
\begin{eqnarray*}
&&\frac{1}{\sqrt{2}}\int_{t_1}^{t_2}\int_{|x|=t}\,\rho\, x^{\beta}\,\frac{x^{\alpha}}{t}\,\partial_{\alpha}u\,\partial_{\beta}u-\rho\,\left(-t\,\partial_{\alpha}u\,\partial^{\alpha}u+x^j\,\frac{x_j}{t}\,\partial_{\alpha}u\,\partial^{\alpha}u\right)\\
&&\quad\quad\quad\quad\quad+\rho \left(-t+\frac{x_jx^j}{t}\right)\frac{|u|^{\dual}}{\dual}+\left(\frac{d}{2}-1\right)\,\rho\,\frac{x^{\alpha}}{t}\,\partial_{\alpha}u\,u\,d\sigma dt\\
&&\quad\quad\quad\quad\quad+\,O\left(\epsilon^{-1}\int_{|x|\leq t_2}|\nabla_{x,t}u|^2+|\nabla_{x,t} u|\,\frac{|u|}{t_2}+|u|^{\dual}(x,t_2)\,dx\right)+\\
&&\quad\quad\quad\quad\quad+\,O\left(\epsilon^{-1}\int_{|x|\leq t_1}|\nabla_{x,t}u|^2+|u|^{\dual}+|\nabla_{x,t} u|\,\frac{|u|}{t_1}(x,t_2)\,dx\right).\\
&&=\frac{1}{\sqrt{2}}\int_{t_1}^{t_2}\int_{|x|=t}\,\rho\, x^{\beta}\,\frac{x^{\alpha}}{t}\,\partial_{\alpha}u\,\partial_{\beta}u+\left(\frac{d}{2}-1\right)\,\rho\,\frac{x^{\alpha}}{t}\,\partial_{\alpha}u\,u\,d\sigma dt+O(\epsilon^{-1}).
\end{eqnarray*}
In the above the $O(\epsilon^{-1})$ denotes a term bounded in magnitude by $C\epsilon^{-1}$. Using the bound on the energy flux (\ref{eq:usebound1}) and estimate (\ref{eq:ineqrho}), we can bound the boundary terms $B$ as 
\begin{eqnarray*}
|B|&\lesssim&\,\int_{t_1}^{t_2}\int_{|x|=t}\,\epsilon^{-1} \left(\frac{x^{\alpha}}{t}\,\partial_{\alpha}u\right)^2\,d\sigma dt+O\left(\epsilon^{-1}\right)\\\\
&&\,\,+\left(\frac{d}{2}-1\right)\left[\int_{t_1}^{t_2}\int_{|x|=t}\,\rho\,t\left(\frac{x^{\alpha}}{t}\,\partial_{\alpha}u\right)^2 d\sigma dt\right]^{\frac{1}{2}}\,\left[\int_{t_1}^{t_2}\int_{|x|=t}\,\rho \,\frac{u^2}{t}\,d\sigma dt\,\right]^{\frac{1}{2}}\\
&\lesssim&\, \epsilon^{-1}+\epsilon^{-\frac{1}{2}}\left(\int_{t_1}^{t_2}\int_{|x|=t}\,\epsilon^{-1}\,t^{-2}\,u^2\,d\sigma dt\,\right)^{\frac{1}{2}}\\
&\lesssim&\, \epsilon^{-1}+\epsilon^{-1}\left(\int_{t_1}^{t_2}\int_{|x|=t}\,|u|^{\dual}\,d\sigma dt\right)^{\frac{1}{\dual}}\,\left(\int_{t_1}^{t_2}\int_{|x|=t}\,t^{-d}\,d\sigma dt\right)^{\frac{1}{d}}\\
&\lesssim&\,\epsilon^{-1}+\epsilon^{-1}\left(\log{\frac{t_2}{t_1}}\right)^{\frac{1}{d}}.
\end{eqnarray*}
The interior terms I are 
\begin{eqnarray*}
I&=&\int_{t_1}^{t_2}\int_{|x|<t}\,-\partial^{\alpha}\left(\rho \,x^{\beta}\right)\,\partial_{\alpha}u\,\partial_{\beta}u+\,\frac{1}{2}\,\partial_{\beta}\left(\rho\,x^{\beta}\right)\partial^{\alpha}u\,\partial_{\alpha}u-\partial_{\beta}\left(\rho\,x^{\beta}\right)\frac{|u|^{\dual}}{\dual}\\
&&\quad\quad\quad-\left(\frac{d}{2}-1\right)\partial^{\alpha}\rho\,\partial_{\alpha}u\,u-\left(\frac{d}{2}-1\right)\rho\,\partial_{\alpha}u\,\partial^{\alpha}u+\left(\frac{d}{2}-1\right)\rho\,|u|^{\dual}\,dxdt,
\end{eqnarray*}
which is
\begin{eqnarray*}
I&=&\int_{t_1}^{t_2}\int_{|x|<t}\,-\partial^{\alpha}\rho\,\partial_{\alpha}u\,x^{\beta}\,\partial_{\beta}u-\rho\,\partial^{\alpha}u\,\partial_{\alpha}u+\frac{1}{2}\,x^{\beta}\,\partial_{\beta}\rho\,\partial^{\alpha}u\,\partial_{\alpha}u\\
&&\quad\quad\quad+\,\frac{d+1}{2}\rho\,\partial_{\alpha}u\,\partial^{\alpha}u-\frac{d+1}{\dual}\rho\,|u|^{\dual}-x^{\beta}\,\partial_{\beta}\rho\,\frac{|u|^{\dual}}{\dual}\\
&&\quad\quad\quad-\left(\frac{d}{2}-1\right)\partial^{\alpha}\rho\,\partial_{\alpha}u\,u-\left(\frac{d}{2}-1\right)\rho\,\partial_{\alpha}u\,\partial^{\alpha}u+\left(\frac{d}{2}-1\right)\rho\,|u|^{\dual}.
\end{eqnarray*}
Recalling that $x^{\beta}\,\partial_{\beta}\rho=-\rho$ and simplifying the above expression, we get that
\begin{equation*}
\int_{t_1}^{t_2}\int_{|x|<t}\,-\partial^{\alpha}\rho\,\partial_{\alpha}u\,x^{\beta}\partial_{\beta}u-\left(\frac{d}{2}-1\right)\partial^{\alpha}\rho\,\partial_{\alpha}u\,u\,dxdt=I=-B=O\left(\epsilon^{-1}\left(\log{\frac{t_2}{t_1}}\right)^{\frac{1}{d}}\right),
\end{equation*}
by the estimates on the boundary terms. Equivalently
\begin{equation*}
\int_{t_1}^{t_2}\int_{|x|<t}\,\partial^{\alpha}\rho\,\partial_{\alpha}u\,x^{\beta}\partial_{\beta}u+\left(\frac{d}{2}-1\right)\partial^{\alpha}\rho\,\partial_{\alpha}u\,u\,dxdt=O\left(\epsilon^{-1}\left(\log{\frac{t_2}{t_1}}\right)^{\frac{1}{d}}\right),
\end{equation*}
that is
\begin{equation}\label{eq:goodid}
\int_{t_1}^{t_2}\int_{|x|<t}\,\partial^{\alpha}\rho\,\partial_{\alpha}u\,\left(x^{\beta}\partial_{\beta}u+\left(\frac{d}{2}-1\right)u\right)\,dxdt=O\left(\epsilon^{-1}\left(\log{\frac{t_2}{t_1}}\right)^{\frac{1}{d}}\right).
\end{equation}
We can rewrite (\ref{eq:goodid}) as
\begin{eqnarray*}
O\left(\epsilon^{-1}\left(\log{\frac{t_2}{t_1}}\right)^{\frac{1}{d}}\right)&=&\int_{t_1}^{t_2}\int_{|x|<t}\,\left(\partial^{\alpha}\rho\,\partial_{\alpha}u+\left(\frac{d}{2}-1\right)\rho^3\,u\right)\,\left(x^{\beta}\partial_{\beta}u+\left(\frac{d}{2}-1\right)u\right)\,dxdt\\
&&-\int_{t_1}^{t_2}\int_{|x|<t}\,\left(\frac{d}{2}-1\right)\rho^3\,u\,\left(x^{\beta}\partial_{\beta}u+\left(\frac{d}{2}-1\right)u\right)\,dxdt.
\end{eqnarray*}
Let us now look at the term $Z:=\int_{t_1}^{t_2}\int_{|x|<t}\,\rho^3\,u\,\left(x^{\beta}\partial_{\beta}u+\left(\frac{d}{2}-1\right)u\right)\,dxdt$. We can calculate
\begin{eqnarray*}
Z&=&\int_{t_1}^{t_2}\int_{|x|<t}\,\rho^3\,x^{\beta}\partial_{\beta}\frac{u^2}{2}+\left(\frac{d}{2}-1\right)\rho^3\,u^2\,dxdt\\
&=&\frac{1}{\sqrt{2}}\,\int_{t_1}^{t_2}\int_{|x|=t}\,\rho^3\,\frac{u^2}{2}\,x^{\beta}\,n_{\beta}\,d\sigma dt+\int_{|x|<t_2}\,\rho^3\,t_2\,\frac{u^2}{2}(x,t_2)\,dx-\int_{|x|<t_1}\,\rho^3\,t_1\,\frac{u^2}{2}(x,t_1)\,dx\\
&&\quad\quad +\int_{t_1}^{t_2}\int_{|x|<t}\,-\partial_{\beta}\left(\rho^3\,x^{\beta}\right)\frac{u^2}{2}+\left(\frac{d}{2}-1\right)\rho^3\,u^2\,dxdt\\
&=&\int_{|x|<t_2}\,\rho^3\,t_2\,\frac{u^2}{2}(x,t_2)\,dx-\int_{|x|<t_1}\,\rho^3\,t_1\,\frac{u^2}{2}(x,t_1)\,dx\\
&&\quad\quad+ \int_{t_1}^{t_2}\int_{|x|<t}\,-3\rho^2\,x^{\beta}\,\partial_{\beta}\rho\,\frac{u^2}{2}-\frac{d+1}{2}\rho^3\,u^2+\left(\frac{d}{2}-1\right)\rho^3\,u^2\,dxdt\\
&=&\int_{|x|<t_2}\,\rho^3\,t_2\,\frac{u^2}{2}(x,t_2)\,dx-\int_{|x|<t_1}\,\rho^3\,t_1\,\frac{u^2}{2}(x,t_1)\,dx\\.
\end{eqnarray*}
Hence we can estimate 
\begin{equation}\label{eq:estimateZ}
|Z|\lesssim \sum_{i=1}^2\int_{|x|<t_i}\epsilon^{-3}\,\frac{u^2}{t_i^2}(x,t_i)\,dx\lesssim \epsilon^{-3}.
\end{equation}
Hence, we have by (\ref{eq:goodid}) and (\ref{eq:estimateZ}) that
\begin{eqnarray*}
&&\int_{t_1}^{t_2}\int_{|x|<t}\,\left(\partial^{\alpha}\rho\,\partial_{\alpha}u+\left(\frac{d}{2}-1\right)\rho^3\,u\right)\,\left(x^{\beta}\partial_{\beta}u+\left(\frac{d}{2}-1\right)u\right)\,dxdt\\
&&=O\left(\epsilon^{-1}\left(\log{\frac{t_2}{t_1}}\right)^{\frac{1}{d}}+\epsilon^{-3}\right).\\
\end{eqnarray*}
Observe that 
\begin{equation}
\partial^{\alpha}\rho\,\partial_{\alpha}u=\rho^3\,\left(x^{\alpha}\partial_{\alpha}u\right)+\epsilon^2\,\rho^3\,t\partial_tu.
\end{equation}
Thus
\begin{eqnarray*}
&&\int_{t_1}^{t_2}\int_{|x|<t}\,\left(\partial^{\alpha}\rho\,\partial_{\alpha}u+\left(\frac{d}{2}-1\right)\rho^3\,u\right)\,\left(x^{\beta}\partial_{\beta}u+\left(\frac{d}{2}-1\right)u\right)\,dxdt\\
&&=\int_{t_1}^{t_2}\int_{|x|<t}\,\rho^3\,\left(x^{\beta}\partial_{\beta}u+\left(\frac{d}{2}-1\right)u\right)^2\,dxdt\\
&&\quad\quad\quad\quad+\,\epsilon^2\,\int_{t_1}^{t_2}\int_{|x|<t}\,\rho^3\,t\,\partial_tu\,\left(x^{\beta}\partial_{\beta}u+\left(\frac{d}{2}-1\right)u\right)\,dxdt\\
&&\ge \frac{1}{2}\,\int_{t_1}^{t_2}\int_{|x|<t}\,\rho^3\,\left(x^{\beta}\partial_{\beta}u+\left(\frac{d}{2}-1\right)u\right)^2\,dxdt\\
&&\quad\quad\quad\quad\quad - \,4\,\epsilon^4\,\int_{t_1}^{t_2}\int_{|x|<t}\,\rho^3\,t^2\,|\partial_tu|^2\,dxdt\\
&&=\frac{1}{2}\,\int_{t_1}^{t_2}\int_{|x|<t}\,\rho^3\,\left(x^{\beta}\partial_{\beta}u+\left(\frac{d}{2}-1\right)u\right)^2\,dxdt\\
&&\quad\quad\quad\quad\quad+O\left(\epsilon \int_{t_1}^{t_2}\int_{|x|<t}\,|\partial_tu|^2\,dx\frac{dt}{t}\right)\\
&&=\frac{1}{2}\,\int_{t_1}^{t_2}\int_{|x|<t}\,\rho^3\,\left(x^{\beta}\partial_{\beta}u+\left(\frac{d}{2}-1\right)u\right)^2\,dxdt+O\left(\epsilon\,\log{\frac{t_2}{t_1}}\right).
\end{eqnarray*}
In summary, we get that 
\begin{eqnarray*}
&&\int_{t_1}^{t_2}\int_{|x|<t}\,\left(\frac{x^{\beta}}{t}\partial_{\beta}u+\left(\frac{d}{2}-1\right)\frac{u}{t}\right)^2\,dx\frac{dt}{t}\\
&&\lesssim \int_{t_1}^{t_2}\int_{|x|<t}\,\rho^3\,\left(x^{\beta}\partial_{\beta}u+\left(\frac{d}{2}-1\right)u\right)^2\,dxdt\\
&&\lesssim\,\epsilon\,\log{\frac{t_2}{t_1}}+\epsilon^{-1}\left(\log{\frac{t_2}{t_1}}\right)^{\frac{1}{d}}+\epsilon^{-3},
\end{eqnarray*}
for any $\epsilon\in(0,1)$. Fix $\epsilon=\left(\log{\frac{t_2}{t_1}}\right)^{-\frac{1}{d+1}}$, we then obtain that
\begin{equation}\label{eq:Morawetz}
\int_{t_1}^{t_2}\int_{|x|<t}\,\left(\partial_tu+\frac{x}{t}\cdot\nabla u+\left(\frac{d}{2}-1\right)\frac{u}{t}\right)^2\,dx\frac{dt}{t}\leq C \left(\log{\frac{t_2}{t_1}}\right)^{\frac{d}{d+1}}.
\end{equation}
The lemma is proved.\\

\hspace{.6cm}This is an extremely useful estimate. Since the right hand side grows slower than $\log{\frac{t_2}{t_1}}$, (\ref{eq:Morawetz}) forces $\int_{|x|<t}\,\left(\partial_tu+\frac{x}{t}\cdot\nabla u+\left(\frac{d}{2}-1\right)\frac{u}{t}\right)^2\,dx$ to vanish on average, which plays an essential role in our analysis below. \\

\end{section}

\begin{section}{Applications of the Morawetz inequality}
The Morawetz estimate implies the following vanishing condition.
\begin{lemma}\label{lm:keyvanishingMorawetz}
Let $u$ be as in the last section. Then there exist $\mu_j\downarrow 0$, $t_j\in (\frac{4}{3}\mu_j,\,\frac{13}{9}\mu_j)$ and $t_j'\in (\frac{14}{9}\mu_j,\,\frac{5}{3}\mu_j)$, such that
$$\frac{1}{\mu_j}\int_{\mu_j}^{2\mu_j}\int_{|x|<t}\,\left(\partial_tu+\frac{x}{t}\cdot\nabla u+\left(\frac{d}{2}-1\right)\frac{u}{t}\right)^2\,dx \,dt\to 0,\,\,{\rm as}\,\,\mu_j\to 0+,$$
and that
\begin{eqnarray*}
&&\sup_{0<\tau<\frac{t_j}{16}}\frac{1}{\tau}\int_{|t_j-t|<\tau}\,\int_{|x|<4t}\,\left(\partial_tu+\frac{x}{t}\cdot\nabla u+\left(\frac{d}{2}-1\right)\frac{u}{t}\right)^2(x,t)\,dx\,dt\\
&&\quad\quad+\sup_{0<\tau<\frac{t_j'}{16}}\frac{1}{\tau}\int_{|t_j'-t|<\tau}\,\int_{|x|<4t}\,\left(\partial_tu+\frac{x}{t}\cdot\nabla u+\left(\frac{d}{2}-1\right)\frac{u}{t}\right)^2(x,t)\,dx\,dt\to 0,\\
&&\,\,{\rm as}\,\,j\to \infty.
\end{eqnarray*}
\end{lemma}

\smallskip
\noindent
{\it Proof.} Recall the following bound for Type II solution $\OR{u}(t)$ in $R^d\times (0,1]$ with $0\in\mathcal{S}$, proved in the last section:
\begin{equation}\label{eq:Morawetz2}
\int_{t_1}^{t_2}\int_{|x|<t}\,\left(\partial_tu+\frac{x}{t}\cdot\nabla u+\left(\frac{d}{2}-1\right)\frac{u}{t}\right)^2\,dx\frac{dt}{t}\leq C \left(\log{\frac{t_2}{t_1}}\right)^{\frac{d}{d+1}}.
\end{equation}
For each large natural number $J$, apply inequality (\ref{eq:Morawetz2}) to $t_1=4^{-J}$ and $t_2=2^{-J}$,  we get that
\begin{equation*}
\sum_{j=0}^{J-1}\,\int_{2^j\,4^{-J}}^{2^{j+1}\,4^{-J}}\int_{|x|<t}\,\left(\partial_tu+\frac{x}{t}\cdot\nabla u+\left(\frac{d}{2}-1\right)\frac{u}{t}\right)^2\,dx\,\frac{dt}{t}\leq C J^{\frac{d}{d+1}}.
\end{equation*}
Hence, there exists $0\leq j\leq J-1$, such that 
\begin{equation*}
\int_{2^j\,4^{-J}}^{2^{j+1}\,4^{-J}}\int_{|x|<t}\,\left(\partial_tu+\frac{x}{t}\cdot\nabla u+\left(\frac{d}{2}-1\right)\frac{u}{t}\right)^2\,dx\,\frac{dt}{t}\leq C J^{-\frac{1}{d+1}}.
\end{equation*}
For a decreasing subsequence $\mu_j$ of $2^j\,4^{-J}$, we get that
\begin{equation}
\frac{1}{\mu_j}\int_{\mu_j}^{2\mu_j}\int_{|x|<t}\,\left(\partial_tu+\frac{x}{t}\cdot\nabla u+\left(\frac{d}{2}-1\right)\frac{u}{t}\right)^2\,dx \,dt\to 0,\,\,{\rm as}\,\,\mu_j\to 0+.
\end{equation}
Since $u(x,t)=v(x,t)$ for $(x,t)\in\{(x,t):\,x\in B_{\delta},\,|x|>t,\,t\in(0,\delta)\}$ and $v$ is regular, we also have
\begin{equation}
\lim_{\mu_j\to 0}\,\frac{1}{\mu_j}\int_{\mu_j}^{2\mu_j}\int_{|x|<4t}\,\left(\partial_tu+\frac{x}{t}\cdot\nabla u+\left(\frac{d}{2}-1\right)\frac{u}{t}\right)^2\,dx \,dt=0.
\end{equation}
Let 
\begin{equation}
g(t)=\int_{|x|<4t}\,\left(\partial_tu+\frac{x}{t}\cdot\nabla u+\left(\frac{d}{2}-1\right)\frac{u}{t}\right)^2(x,t)\,dx,
\end{equation}
then $\frac{1}{\mu_j}\int_{\mu_j}^{2\mu_j}\,g(t)\,dt\to 0+$ as $j\to\infty$. Denote $M(g\chi_{(\mu_j,\,2\mu_j)})$ as the Hardy-Littlewood maximal function of $g\chi_{(\mu_j,\,2\mu_j)}$. Passing to a subsequence, we can assume that 
\begin{equation*}
\frac{1}{\mu_j}\int_{\mu_j}^{2\mu_j}\,g(t)\,dt\leq 4^{-j}.
\end{equation*}
Then using the fact that the Hardy-Littlewood maximal operator $M$ is bounded from $L^1$ to $L^{1,\infty}$, we get that
\begin{equation*}
\left|\{t\in (\mu_j,\,2\mu_j):\,M(g\chi_{(\mu_j,\,2\mu_j)})(t)>2^{-j}\}\right|\leq C 2^{-j}\mu_j.
\end{equation*}
Thus we can find sequences $t_j$ and $t_j'$ such that $t_j\in (\frac{4}{3}\mu_j,\,\frac{13}{9}\mu_j)$ and $t_j'\in (\frac{14}{9}\mu_j,\,\frac{5}{3}\mu_j)$, such that $M(g\chi_{(\mu_j,\,2\mu_j)})(t_j)\to 0+$ and $M(g\chi_{(\mu_j,\,2\mu_j)})(t_j')\to 0+$ as $j\to\infty$. Consequently along sequences $0<t_j<t_j'$, we have 
\begin{equation}
\sup_{0<\tau<\frac{t_j}{16}}\frac{1}{\tau}\int_{|t_j-t|<\tau}\,g(t)\,dt+\sup_{0<\tau<\frac{t_j'}{16}}\frac{1}{\tau}\int_{|t_j'-t|<\tau}\,g(t)\,dt\to 0,\,\,{\rm as}\,\,j\to \infty.
\end{equation}
That is 
\begin{eqnarray*}
&&\sup_{0<\tau<\frac{t_j}{16}}\frac{1}{\tau}\int_{|t_j-t|<\tau}\,\int_{|x|<4t}\,\left(\partial_tu+\frac{x}{t}\cdot\nabla u+\left(\frac{d}{2}-1\right)\frac{u}{t}\right)^2(x,t)\,dx\,dt\\
&&\quad\quad+\sup_{0<\tau<\frac{t_j'}{16}}\frac{1}{\tau}\int_{|t_j'-t|<\tau}\,\int_{|x|<4t}\,\left(\partial_tu+\frac{x}{t}\cdot\nabla u+\left(\frac{d}{2}-1\right)\frac{u}{t}\right)^2(x,t)\,dx\,dt\to 0,\\
&&\,\,{\rm as}\,\,j\to \infty.
\end{eqnarray*}
The lemma is proved.\\

\hspace{.6cm}The reason why we need two sequences $t_j$, $t_j'$ approaching zero will become clear below.\\

\end{section}

\begin{section}{Characterization of solutions along a time sequence with the vanishing condition}

In this section, we use the asymptotic vanishing of the quantity 
\begin{equation}
\int_{|x|<4t}\,\left(\partial_tu+\frac{x}{t}\cdot\nabla u+\left(\frac{d}{2}-1\right)\frac{u}{t}\right)^2(x,t)\,dx,
\end{equation}
to obtain a preliminary decomposition along a sequence of times in the spirit of decomposition (\ref{eq:maindecompositionhaha}), albeit with a remainder term that vanishes only in $L^{\dual}$.  Suppose that $\OR{u}\in C((0,1],\HL)$, with $u\in \Snorm(R^d\times (\epsilon,1])$ for any $\epsilon>0$, is a type II solution to equation (\ref{eq:main}). We assume that $x_{\ast}=0$ is a singular point. Then by Lemma \ref{lm:keyvanishingMorawetz}, there exists a sequence of times $t_n\downarrow 0$, such that
\begin{equation}\label{eq:vanishingcharacterization}
\lim_{n\to\infty}\,\sup_{0<\tau<\frac{t_n}{16}}\frac{1}{\tau}\int_{|t_n-t|<\tau}\,\int_{|x|<4t}\,\left(\partial_tu+\frac{x}{t}\cdot\nabla u+\left(\frac{d}{2}-1\right)\frac{u}{t}\right)^2(x,t)\,dx\,dt = 0.
\end{equation}
Our main goal in this section is to prove the following theorem.
\begin{theorem}\label{th:mainpreliminary}
Let $\OR{u}$, $t_n$ be as above. Let $\OR{v}\in C([0,\,\delta],\HL)$, with $v\in \Snorm(R^d\times [0,\delta])$ for some $\delta>0$, be the regular part of $\OR{u}$. Then passing to a subsequence there exist integer $J_0\ge 0$\,\footnote{If $J_0=0$, then there is no sum of solitons. We shall show in the last section that $J_0\ge 1$.}, $r_0>0$, scales $\lambda_n^j$ with $0<\lambda_n^j\ll t_n$, positions $c_n^j\in R^d$ satisfying $c_n^j\in B_{\beta\, t_n}$ for some $\beta\in(0,1)$, with $\ell_j=\lim\limits_{n\to\infty}\frac{c_n^j}{t_n}$ well defined, and travalling waves $Q_{\ell_j}$, for $1\leq j\leq J_0$, such that
\begin{equation}\label{eq:maindecomposition}
\OR{u}(t_n)=\OR{v}(t_n)+\sum_{j=1}^{J_0}\,\left((\lambda_n^j)^{-\frac{d}{2}+1}\, Q_{\ell_j}\left(\frac{x-c_n^j}{\lambda_n^j},\,0\right),\,(\lambda_n^j)^{-\frac{d}{2}}\, \partial_tQ_{\ell_j}\left(\frac{x-c_n^j}{\lambda_n^j},\,0\right)\right)+(\epsilon_{0n},\epsilon_{1n}),
\end{equation}
where $(\epsilon_{0n},\epsilon_{1n})$ vanishes asymptotically in the following sense:
\begin{equation}
\|\epsilon_{0n}\|_{L^{\dual}(|x|\leq r_0)}\to 0,\quad {\rm as}\,\,n\to\infty;
\end{equation}
if we write $\OR{\epsilon}_n=(\epsilon_{0n},\,\epsilon_{1n})$ as
\begin{equation}
\OR{\epsilon}_n(x)=\left((\lambda^j_n)^{-\frac{d}{2}+1}\widetilde{\epsilon}_{0n}\left(\frac{x-c^j_n}{\lambda^j_n}\right),\,(\lambda^j_n)^{-\frac{d}{2}}\widetilde{\epsilon}_{1n}\left(\frac{x-c^j_n}{\lambda^j_n}\right)\right),
\end{equation}
then $\OR{\widetilde{\epsilon}}_n\rightharpoonup 0$ as $n\to\infty$, for each $j\leq J_0.$
In addition, the parameters $\lambda_n^j,\,c^j_n$ satisfy the pseudo-orthogonality condition
\begin{equation}\label{eq:pseudo}
\frac{\lambda_n^j}{\lambda_n^{j'}}+\frac{\lambda_n^{j'}}{\lambda_n^j}+\frac{\left|c_n^j-c_n^{j'}\right|}{\lambda_n^j}\to\infty,
\end{equation}
as $n\to\infty$, for $1\leq j\neq j'\leq J_0.$
\end{theorem}

\smallskip
\noindent
{\it Proof.} Fix $\phi\in C_c^{\infty}(B_4)$ with $\phi|_{B_3}\equiv 1$. Let 
\begin{equation*}
(u_{0n},\,u_{1n})=\OR{u}(t_n)\phi\left(\frac{x}{t_n}\right).
\end{equation*}
Clearly $(u_{0n},\,u_{1n})$ is a bounded sequence in $\HL$.
Passing to a subsequence if necessary, we can assume that $(u_{0n},\,u_{1n})$ has the profile decomposition 
\begin{eqnarray*}
&&(u_{0n},\,u_{1n})\\
&&\quad=\sum_{j=1}^J\,\left(\left(\lambda^j_n\right)^{-\frac{d}{2}+1}U_j^L\left(\frac{x-c^j_n}{\lambda^j_n}, \,-\frac{t^j_n}{\lambda^j_n}\right),\,\left(\lambda^j_n\right)^{-\frac{d}{2}}\partial_tU_j^L\left(\frac{x-c^j_n}{\lambda^j_n}, \,-\frac{t^j_n}{\lambda^j_n}\right)\right)+(w^J_{0n},\,w^J_{1n}),
\end{eqnarray*}
where the parameters satisfy
\begin{equation}\label{profileorthogonality3}
t^j_n\equiv 0\,\,{\rm for\,\,all\,\,}n,\,\,\,\,{\rm or}\,\,\lim_{n\to\infty}\frac{t^j_n}{\lambda^j_n}\in\{\pm \infty\},
\end{equation}
and for $j\neq j'$
\begin{equation}\label{eq:profileorthogonality4}
\lim_{n\to\infty}\,\left(\frac{\lambda^j_n}{\lambda_n^{j'}}+\frac{\lambda^{j'}_n}{\lambda^j_n}+\frac{\left|c^j_n-c^{j'}_n\right|}{\lambda^j_n}+\frac{\left|t^j_n-t^{j'}_n\right|}{\lambda^j_n}\right)=\infty.
\end{equation}
We note that $(u_{0n},\,u_{1n})$ has the following concentration property
\begin{equation}\label{eq:concentrationproperty}
\lim_{n\to\infty}\,\|(u_{0n},\,u_{1n})\|_{\HL(|x|\ge t_n)}=0.
\end{equation}
By the concentration property (\ref{eq:concentrationproperty}) of $(u_{0n},\,u_{1n})$ and the linear profiles, we can assume that $\lambda^j_n\lesssim t_n$, and if $\lambda^j_n\sim t_n$ along a subsequence, then passing to a subsequence and rescaling the profile, we choose $\lambda^j_n=t_n$. Additionally, if $t^j_n\equiv 0$, we can assume that $\frac{c^j_n}{t_n}$ is bounded, and by passing to a subsequce also that $\ell_j=\lim\limits_{n\to\infty}\frac{c^j_n}{t_n}$ is well defined for such $j$.
We divide the profiles into three cases.\\
\begin{itemize}

\item {\it Case I:} \,\,$t^j_n\equiv 0$, $\lambda^j_n\equiv t_n$. We shall show that the nonlinear profile $U_j$ is a compactly supported self similar solution, which we know must be trivial from \cite{DKMacta};\\

\item {\it Case II:} \,\, $t^j_n\equiv 0$, $\lambda^j_n\ll t_n$. In this case, we shall show that $|\ell_j|<1$ and $\OR{U}_j^L(\cdot,0)=\OR{Q}_{\ell_j}$;\\

\item {\it Case III:}\,\, $\lambda^j_n\ll \left|t^j_n\right|$. We show that these profiles can be absorbed into the residue term. \\

\end{itemize}

Let us firstly consider {\it Case I}. To simplify notations, we assume that the profile is $U^L$. By the support property of $(u_{0n},\,u_{1n})$, we can also assume that the position is $c^j_n\equiv 0$, then we see that ${\rm supp}\,\OR{U}^L\subseteq \overline{B_1}$.  Denote the nonlinear profile as $U$. Assume that $U$ exists in $R^d\times[-T,T]$, with $\|U\|_{\Snorm(R^d\times[-T,T])}\leq M<\infty$. Fix $\epsilon>0$ sufficiently small, determined by the small data Cauchy theory for equation (\ref{eq:main}). By the energy expansion for the profile decomposition, there exists $J(\epsilon)$ such that
\begin{equation}
\sum_{j\ge J(\epsilon)}\|\OR{U}^L_j\|_{\HL}^2<\epsilon^2.
\end{equation}
The profiles with $j\ge J(\epsilon)$ are small, and will be controlled perturbatively by Lemma \ref{lm:nonlinearprofiledecomposition}. For the remaining profiles, passing to a subsequence if necessary, we can assume that $\lim\limits_{n\to\infty}\frac{t^j_n}{t_n}\in\{\pm\infty\}$, or $\lim\limits_{n\to\infty}\frac{t^j_n}{t_n}$ exists and is finite.  We destinguish two categories:\\

\begin{itemize}

\item $j\leq J(\epsilon)$, $\lim\limits_{n\to\infty}\frac{t^j_n}{t_n}\neq 0$. Then for some $\gamma>1$ sufficiently large, we have
\begin{eqnarray*}
&&\left\|\left(\left(\lambda^j_n\right)^{-\frac{d}{2}+1}U_j^L\left(\frac{x-c^j_n}{\lambda^j_n}, \,\frac{t-t^j_n}{\lambda^j_n}\right),\right.\right.\\
&&\quad\quad\quad\quad\quad\quad\,\left.\left.\left(\lambda^j_n\right)^{-\frac{d}{2}}\partial_tU_j^L\left(\frac{x-c^j_n}{\lambda^j_n}, \,\frac{t-t^j_n}{\lambda^j_n}\right)\right)\right\|_{\Snorm\left(R^d\times(-\frac{t_n}{\gamma},\,\frac{t_n}{\gamma})\right)}\to 0,
\end{eqnarray*}
as $n\to\infty$. These profiles have asymptotically vanishing interaction with the profile $U$, and can be controlled perturbatively at least for a short time that is comparable to $t_n$;\\

\item $t^j_n=o(t_n)$, or, $t^j_n\equiv 0$ and $\lambda^j_n\ll t_n$. By the concentration property of $(u_{0n},\,u_{1n})$, passing to a subsequence, we can assume that $\lim\limits_{n\to\infty}\frac{c^j_n}{t_n}=\ell_j$ is well defined. These profiles may have nontrivial interaction with the first profile, and will be removed. Denote the set of such $j$ as $\mathcal{J}_1$.\\

\end{itemize}
Recall that
\begin{equation*}
\OR{U}^L_{jn}(\cdot,0):=\left(\left(\lambda^j_n\right)^{-\frac{d}{2}+1}U_j^L\left(\frac{x-c^j_n}{\lambda^j_n}, \,-\frac{t^j_n}{\lambda^j_n}\right),\,\left(\lambda^j_n\right)^{-\frac{d}{2}}\partial_tU_j^L\left(\frac{x-c^j_n}{\lambda^j_n}, \,-\frac{t^j_n}{\lambda^j_n}\right)\right).
\end{equation*}
For any $\tau>0$ and $j\in \mathcal{J}_1$, by Lemma \ref{lm:concentrationoffreeradiation}, one can verify that
\begin{equation}
\left\|\OR{U}^L_{jn}(0)\right\|_{\HL\left(R^d\backslash\,B_{\tau t_n}(\ell_jt_n)\right)}+\|U^L_{jn}(0)\|_{L^{\dual}\left(R^d\backslash\,B_{\tau t_n}(\ell_jt_n)\right)}\longrightarrow 0,\,\,\,\,{\rm as}\,\,n\to\infty.\label{eq:closenesslocal}
\end{equation}
Define 
\begin{eqnarray*}
(v_{0n},\,v_{1n})&=&\left(t_n^{-\frac{d}{2}+1}U_0\left(\frac{x}{t_n}\right),\,t_n^{-\frac{d}{2}}U_1\left(\frac{x}{t_n}\right)\right)\\
&&\quad+\sum_{j\leq J,\,j\not\in \mathcal{J}_1}\left(\left(\lambda^j_n\right)^{-\frac{d}{2}+1}U_j^L\left(\frac{x-c^j_n}{\lambda^j_n}, \,-\frac{t^j_n}{\lambda^j_n}\right),\,\left(\lambda^j_n\right)^{-\frac{d}{2}}\partial_tU_j^L\left(\frac{x-c^j_n}{\lambda^j_n}, \,-\frac{t^j_n}{\lambda^j_n}\right)\right)\\
&&\\
&&\quad\quad+\,(w^J_{0n},\,w^J_{1n}).
\end{eqnarray*}
Consider the rescaled sequences
\begin{eqnarray*}
&&(\widetilde{u}_{0n},\,\widetilde{u}_{1n})=\left(t_n^{\frac{d}{2}-1}\,u_{0n}(t_n\,x),\,t_n^{\frac{d}{2}}\,u_{1n}(t_n\,x)\right);\\
&&(\widetilde{v}_{0n},\,\widetilde{v}_{1n})=\left(t_n^{\frac{d}{2}-1}\,v_{0n}(t_n\,x),\,t_n^{\frac{d}{2}}\,v_{1n}(t_n\,x)\right).
\end{eqnarray*}
Hence
\begin{eqnarray*}
\left(\widetilde{v}_{0n},\,\widetilde{v}_{1n}\right)&=&(U_0,\,U_1)+\,\left(t_n^{\frac{d}{2}-1}\,w^J_{0n}(t_nx),\,t_n^{\frac{d}{2}}\,w^J_{1n}(t_nx)\right)\\
&&\quad\quad\quad+\,\sum_{j\leq J,\,j\not\in \mathcal{J}_1}\left(\left(t_n^{-1}\lambda^j_n\right)^{-\frac{d}{2}+1}U_j^L\left(\frac{x-t_n^{-1}c^j_n}{t_n^{-1}\lambda^j_n}, \,-\frac{t_n^{-1}t^j_n}{t_n^{-1}\lambda^j_n}\right),\right.\\
&&\quad\quad\quad\quad\quad\quad\quad\quad\quad\quad\left.\left(t_n^{-1}\lambda^j_n\right)^{-\frac{d}{2}}\partial_tU_j^L\left(\frac{x-t_n^{-1}c^j_n}{t_n^{-1}\lambda^j_n}, \,-\frac{t_n^{-1}t^j_n}{t_n^{-1}\lambda^j_n}\right)\right).
\end{eqnarray*}
Let $\widetilde{u}_n$ and $\widetilde{v}_n$ be the solutions to equation (\ref{eq:main}) with intial data $(\widetilde{u}_{0n},\,\widetilde{u}_{1n})$ and $(\widetilde{v}_{0n},\,\widetilde{v}_{1n})$ respectively. 

  \hspace{.6cm}    By (\ref{eq:closenesslocal}), and a rescaling, we see that
\begin{equation}\label{eq:rescaledcloseness}
\left\|\left(\widetilde{u}_{0n},\,\widetilde{u}_{1n}\right)-\left(\widetilde{v}_{0n},\,\widetilde{v}_{1n}\right)\right\|_{\HL\left(R^d\backslash\,\bigcup_{j\in\mathcal{J}_1}B_{\tau }(\ell_j)\right)}\to 0,
\end{equation}
as $n\to\infty$, for any $\tau>0$.

    \hspace{.6cm}For $T_1$ sufficiently small, by Lemma \ref{lm:nonlinearprofiledecomposition}, $\OR{\widetilde{v}}_n$ exists in $R^d\times [-T_1,\,T_1]$, with bound
\begin{equation*}
\|\widetilde{v}_n\|_{\Snorm\left(R^d\times[-T_1,T_1]\right)}\leq M_1<\infty.
\end{equation*}
For any $R>0$, $\sigma>0$ and any $y$ with 
\begin{equation*}
\inf\,\{|y-\ell_j|,\,j\in \mathcal{J}_1\}\ge R+\sigma,
\end{equation*}
(\ref{eq:rescaledcloseness}) implies that
\begin{equation*}
\|(\widetilde{u}_{0n},\,\widetilde{u}_{1n})-(\widetilde{v}_{0n},\,\widetilde{v}_{1n})\|_{\HL(B_{R+\frac{\sigma}{2}}(y))}+\|\widetilde{u}_{0n}-\widetilde{v}_{0n}\|_{L^{\dual}(B_{R+\frac{\sigma}{2}}(y))}\to 0,\,\,{\rm as}\,\,n\to\infty.
\end{equation*}
Hence by Lemma \ref{lm:localinspaceapproximation}, for $R_1=\min\left\{T_1,\,R\right\}$, we get that
\begin{eqnarray}
&&\sup_{t\in[-R_1,\,R_1]}\,\|\OR{\widetilde{v}}_n(t)-\OR{\widetilde{u}}_n(t)\|_{\HL(|x-y|\leq R-|t|)}\nonumber\\
&&\nonumber\\
&&\quad\quad\quad\quad\quad+\,\|\widetilde{v}_n-\widetilde{u}_n\|_{\Snorm\left(\left\{(x,t):\,|x|\leq R-|t|,\,|t|\leq R_1\right\}\right)}\to 0+.\label{eq:localcloseness}
\end{eqnarray}
On the other hand, by Lemma \ref{lm:nonlinearprofiledecomposition}, the solution $\OR{\widetilde{v}}_n$ admits the following expansion for $|t|\leq R_1$
\begin{eqnarray*}
\OR{\tilde{v}}_n&=&\OR{U}(x,t)+\left(t_n^{\frac{d}{2}-1}\,w^J_{n}(t_nx,\,t_nt),\,t_n^{\frac{d}{2}}\,\partial_tw^J_{n}(t_nx,\,t_nt)\right)+\OR{r}^J_n\\
&&\quad\quad\quad\quad+\sum_{j\not\in \mathcal{J}_1}\left(\left(t_n^{-1}\lambda^j_n\right)^{-\frac{d}{2}+1}U_j\left(\frac{x-t_n^{-1}c^j_n}{t_n^{-1}\lambda^j_n}, \,\frac{t-t_n^{-1}t^j_n}{t_n^{-1}\lambda^j_n}\right),\right.\\
&&\quad\quad\quad\quad\quad\quad\quad\quad\quad\left.\left(t_n^{-1}\lambda^j_n\right)^{-\frac{d}{2}}\partial_tU_j\left(\frac{x-t_n^{-1}c^j_n}{t_n^{-1}\lambda^j_n}, \,\frac{t-t_n^{-1}t^j_n}{t_n^{-1}\lambda^j_n}\right)\right),
\end{eqnarray*}
where $\lim\limits_{J\to\infty}\limsup\limits_{n\to\infty}\,\sup_{|t|\leq R_1}\,\|\OR{r}^J_n(t)\|_{\HL}=0$.
By pseudo-orthogonality of parameters (\ref{eq:profileorthogonality4}) and the approximation (\ref{eq:localcloseness}), we get that
\begin{equation}\label{11111}
 \OR{\widetilde{u}}_n\rightharpoonup \OR{U}\,\,\,\,{\rm in}\,\, \{(x,t):\,|x-y|\leq R-|t|,\,|t|\leq R_1\}.
\end{equation} 
By finite speed of propagation and rescaling, we have that
\begin{equation}
\widetilde{u}_n(x,t)=t_n^{\frac{d}{2}-1}u(t_nx,\,t_n+t_nt),
\end{equation}
for $|x|<2$ and $|t|<1$.
After rescaling, (\ref{11111}) and (\ref{eq:vanishingcharacterization}) with $\tau=\frac{t_n}{20}$ imply that
\begin{equation*}
\int_{\left\{(x,t):\,|x-y|\leq R-|t|,\,|t|\leq R_1\right\}}\left(\partial_tU+\frac{x}{t+1}\cdot\nabla U+\left(\frac{d}{2}-1\right) \frac{U}{t+1}\right)^2(x,t)\,dxdt=0.
\end{equation*}
Hence 
\begin{equation*}
\partial_tU+\frac{x}{t+1}\cdot\nabla U+\left(\frac{d}{2}-1\right) \frac{U}{t+1}\equiv 0,\,\,{\rm in}\,\,\{(x,t):\,|x|<2,\,|x-y|\leq R-|t|,\,|t|\leq R_1\}.
\end{equation*}
Moving $y$ around and thrinking $\sigma$ to $0$, we can conclude that in fact
\begin{equation}\label{eq:firstorderselfsimilar}
\partial_tU+\frac{x}{t+1}\cdot\nabla U+\left(\frac{d}{2}-1\right) \frac{U}{t+1}\equiv 0,
\end{equation}
in $\left\{(x,t):\,|t|<{\rm dist}\left(x,\,\{\ell_j,\,j\in \mathcal{J}_1\}\right),\,\,|x|<2,\,\,|t|<R_1\right\}$. Equation (\ref{eq:firstorderselfsimilar}) implies that for some $\Psi$
\begin{equation}\label{eq:selfsimilarrepresentation}
U(x,t)=(t+1)^{-\frac{d}{2}+1}\Psi\left(\frac{x}{t+1}\right),
\end{equation}
in $\left\{(x,t):\,|t|<\frac{1}{4}\,{\rm dist}\left(x,\,\{\ell_j,\,j\in \mathcal{J}_1\}\right),\,|t|<R_1,\,|x|<2\right\}.$ Since ${\rm supp}\,\OR{U}(0)\subseteq \overline{B_1}$, we conclude that 
\begin{equation}\label{eq:equalinitialdata}
\OR{U}(0)=\left(\Psi,\,-x\cdot\nabla \Psi-\left(\frac{d}{2}-1\right)\Psi\right)\,\,\,{\rm in}\,\,R^d.
\end{equation}
Since $\OR{U}$ is supported in $\overline{B_1}$ and satisfies equation (\ref{eq:main}) in $R^d\times [-T,\,T]$, we get from the above that ${\rm supp}\,\Psi\subseteq \overline{B_1}$ and
\begin{equation}\label{eq:degenerateelliptic1}
-\Delta \Psi+y\cdot\nabla (y\cdot\nabla \Psi)+(d-1)y\cdot\nabla \Psi+\frac{d(d-2)}{4}\Psi=|\Psi|^{\dual-2}\Psi,
\end{equation} 
in $R^d\backslash\{\ell_j,\,j\in \mathcal{J}_1\}$. By the regularity condition that $\Psi\in \dot{H}^1$, which follows directly from $\OR{U}\in\HL$ and the representation (\ref{eq:selfsimilarrepresentation}), we can conclude that equation (\ref{eq:degenerateelliptic1}) is in fact satisfied ``across" the points $\ell_j$, i.e., $\ell_j$ are removable singularities. Thus equation (\ref{eq:degenerateelliptic1}) holds in $R^d$. By standard elliptic regularity theory, we see that $\Psi\in C^{2}(B_1)$. Hence  
$$\widetilde{U}(x,t):=(t+1)^{-\frac{d}{2}+1}\Psi\left(\frac{x}{t+1}\right)$$
is a classical solution to equation (\ref{eq:main}) for $|x|<1-|t|$. By (\ref{eq:equalinitialdata}), we see that $\OR{U}(0)=\OR{\widetilde{U}}(0)$. Hence, by finite speed of propagation, $\widetilde{U}\equiv U$ for $|x|<1-|t|$. Consequently $\widetilde{U}\in \Snorm(|x|<1-|t|)$. By the support property of $\Psi$, thus $\widetilde{U}\in\Snorm(R^d\times [-1/2,1/2])$ and is a solution to equation (\ref{eq:main}) with the same initial data as that of $U$. Therefore $U\equiv \widetilde{U}$ is a compactly supported finite energy self similar solution. By \cite{DKMacta}, we see that $U$ must be trivial. \\

\smallskip
\noindent
{\it Remark.} Since we have an exact self similar solution, in principle, we might be able to use simpler arguments than those in \cite{DKMacta} to prove that $U$ must be trivial. \\


\smallskip
\noindent
{\it Case II:} Now we consider Case II, $t^j_n\equiv 0$, $\lambda^j_n\ll t_n$. To simplify notations, we assume that the profile is $U_1^L$ with parameters $\lambda^1_n$, $c^1_n$ satisfying $\lim\limits_{n\to\infty}\frac{c^1_n}{t^1_n}=\ell_1$. 
Suppose that the nonlinear profile $U_1$ associated with $U_1^L$, $\lambda^1_n$ exists in $R^d\times [-T,\,T]$ with 
\begin{equation*}
\|U_1\|_{\Snorm(R^d\times[-T,T])}\leq M<\infty.
\end{equation*}
The idea in the characterization of $U_1$ is similar to Case I. We shall still remove the profiles which contain more energy than the threshold energy provided by small data theory and which have nontrivial interactions with $U_1$, (these profiles are only finitely many), and then use perturbative arguments to deal with other small profiles. After a rescaling, we can still pass to limit in the region with several lightcones removed. Using the vanishing condition (\ref{eq:vanishingcharacterization}), we then obtain a first order equation for $U_1$ which will enable us to classify $U_1$. The difference with Case I here is that due to the fact that $\lambda^1_n\ll t_n$, the first order equation we obtain in the end is different from the self similar case, and as a consequence $U_1$ has to be a travelling wave, instead of a self similar solution. 

\hspace{.6cm}Fix $\epsilon=\epsilon(d)>0$ sufficiently small, determined by the small data Cauchy theory for equation (\ref{eq:main}), we can find large natural number $J(\epsilon)$, such that
\begin{equation*}
\sum_{j\ge J(\epsilon)}\|\OR{U}_j^L\|_{\HL}^2\leq \epsilon^2.
\end{equation*}
These profiles can be controlled perturbatively. For the remaining profiles, we distinguish several categories.\\

\begin{itemize}

\item $1<j<J(\epsilon)$, the profile satisfies the properties that $\lambda^j_n\ll \left|t^j_n\right|$, and (passing to a subsequence if necessary) that $\left|t^j_n\right|\gtrsim \lambda^1_n$ for all $n$. These profiles have negligible interference with the profile $U_1^L$, $\lambda^1_n$ at least for a short time (comparable to $\lambda^1_n$), and can be controlled perturbatively. More precisely, we have
\begin{equation}\label{eq:largeprofileok}
\lim_{n\to\infty}\,\left\|(\lambda^j_n)^{-\frac{d}{2}+1}U^L_j\left(\frac{x-c^j_n}{\lambda^j_n},\,\frac{t-t^j_n}{\lambda^j_n}\right)\right\|_{\Snorm\left(R^d\times [-T_1\lambda^1_n,\,T_1\lambda^1_n]\right)}\to 0,
\end{equation}
as $n\to\infty$, if $T_1>0$ is sufficiently small; \\

\item $1<j<J(\epsilon)$ and the profile satisfies the property that $\lambda^j_n\ll \left|t^j_n\right|$, $t^j_n=o(\lambda^1_n)$ as $n\to\infty$. Denote the set of these $j$ as $\mathcal{J}_1$. These profiles may have nontrivial interaction with the profile $U_1^L$, $\lambda^1_n$, and have to be removed;\\

\item $1<j<J(\epsilon)$, and the profile satisfies the property that $t^j_n\equiv 0$,  and passing to a subsequence if necessary that $\lambda^j_n\ll \lambda^1_n$. Denote the set of such $j$ as $\mathcal{J}_2$. These profiles also have to be removed.\\

\item $1<j<J(\epsilon)$, $t^j_n\equiv 0$, and passing to a subsequence if necessary that $\lambda^j_n\sim \lambda^1_n$, $\frac{\left|c^j_n-c^1_n\right|}{\lambda^1_n}\to\infty$. These profiles have no interaction with the first profile locally in space in the limit.\\

\item $1<j<J(\epsilon)$, and $t^j_n\equiv 0$, $\lambda^1_n=o(\lambda^j_n)$. These profiles have asymtotically vanishing interaction with the first profile for a time interval $\sim \lambda^1_n$, and can be controlled perturbatively.\\

\end{itemize}

Define
\begin{eqnarray*}
(v_{0n},\,v_{1n})&=&\left((\lambda^1_n)^{-\frac{d}{2}+1}U^L_1\left(\frac{x-c^1_n}{\lambda^1_n},\,0\right),\,(\lambda^1_n)^{-\frac{d}{2}}\partial_tU_1^L\left(\frac{x-c^1_n}{\lambda^1_n},\,0\right)\right)\\
&&\quad+\sum_{1<j\leq J,\,\,\,j\not\in \mathcal{J}_1\cup \mathcal{J}_2}\left(\left(\lambda^j_n\right)^{-\frac{d}{2}+1}U_j^L\left(\frac{x-c^j_n}{\lambda^j_n}, \,-\frac{t^j_n}{\lambda^j_n}\right),\,\left(\lambda^j_n\right)^{-\frac{d}{2}}\partial_tU_j^L\left(\frac{x-c^j_n}{\lambda^j_n}, \,-\frac{t^j_n}{\lambda^j_n}\right)\right)\\
&&\\
&&\quad\quad\quad\quad+\,(w^J_{0n},\,w^J_{1n}).
\end{eqnarray*}
We can check that
\begin{eqnarray}
&&\|(v_{0n},\,v_{1n})-(u_{0n},\,u_{1n})\|_{\HL\left(B_{M\lambda^1_n}(c^1_n)\backslash\, \bigcup_{j\in \mathcal{J}_1\cup \mathcal{J}_2} B_{\epsilon \lambda^1_n}(c^j_n)\right)}\nonumber\\
&&\quad\quad\quad\quad+\|v_{0n}-u_{0n}\|_{L^{\dual}\left(B_{M\lambda^1_n}(c^1_n)\backslash\, \bigcup_{j\in \mathcal{J}_1\cup \mathcal{J}_2} B_{\epsilon \lambda^1_n}(c^j_n)\right)}\,\to\, 0,\label{eq:errorvanishesinitial}
\end{eqnarray}
as $n\to\infty$, for any $M>1$ and $\epsilon>0$. \\
Consider the rescaled sequences
\begin{eqnarray*}
&&(\widetilde{u}_{0n},\,\widetilde{u}_{1n})=\left(t_n^{\frac{d}{2}-1}u_{0n}(t_n\,x),\,t_n^{\frac{d}{2}}u_{1n}(t_n\,x)\right);\\
&&(\widetilde{v}_{0n},\,\widetilde{v}_{1n})=\left(t_n^{\frac{d}{2}-1}v_{0n}(t_n\,x),\,t_n^{\frac{d}{2}}v_{1n}(t_n\,x)\right).
\end{eqnarray*}
Let $\widetilde{u}_n$, $\widetilde{v}_n$ be the solution to equation (\ref{eq:main}) with initial data $(\widetilde{u}_{0n},\,\widetilde{u}_{1n})$ and $(\widetilde{v}_{0n},\,\widetilde{v}_{1n})$, respectively. By the principle of finite speed of propagation and rescaling, we get that
\begin{equation}
\widetilde{u}_n(x,t)=t_n^{\frac{d}{2}-1}\,u(t_nx,\,t_n(t+1)),\,\,{\rm for}\,\,|x|<2,\,\,|t|<\frac{1}{2}.
\end{equation}
By Lemma \ref{lm:nonlinearprofiledecomposition}, $\widetilde{v}_n$ admits the following decomposition for $|t|<T_1\lambda^1_n$
\begin{eqnarray*}
\widetilde{v}_n&=&\left((t_n^{-1}\lambda^1_n)^{-\frac{d}{2}+1}U_1\left(\frac{x-t_n^{-1}c^1_n}{t_n^{-1}\lambda^1_n},\,\frac{t}{t_n^{-1}\lambda^1_n}\right),\,(t_n^{-1}\lambda^1_n)^{-\frac{d}{2}}\partial_tU_1\left(\frac{x-t_n^{-1}c^1_n}{t_n^{-1}\lambda^1_n},\,\frac{t}{t_n^{-1}\lambda^1_n}\right)\right)\\
&&\quad\quad\quad+\sum_{1<j\leq J,\,\,j\not\in \mathcal{J}_1\cup \mathcal{J}_2}\left(\left(t_n^{-1}\lambda^j_n\right)^{-\frac{d}{2}+1}U_j\left(\frac{x-t_n^{-1}c^j_n}{t_n^{-1}\lambda^j_n}, \,\frac{t-t_n^{-1}t^j_n}{t_n^{-1}\lambda^j_n}\right),
\right.\\
&&\quad\quad\quad\quad\quad\quad\quad\quad\quad\quad\quad\quad\,\left.\left(t_n^{-1}\lambda^j_n\right)^{-\frac{d}{2}}\partial_tU_j\left(\frac{x-t_n^{-1}c^j_n}{t_n^{-1}\lambda^j_n}, \,\frac{t-t_n^{-1}t^j_n}{t_n^{-1}\lambda^j_n}\right)\right)\\
&&\\
&&\quad\quad\quad\quad\quad+\left(t_n^{\frac{d}{2}-1}w_n^J(t_nx,t_nt),\,t_n^{\frac{d}{2}}\partial_tw^J_n(t_nx,t_nt)\right)+\OR{r}^J_n,
\end{eqnarray*}
where $\OR{r}^J_n$ vanishes asymptotically:
\begin{equation}
\lim_{J\to\infty}\limsup_{n\to\infty}\,\sup_{|t|<T_1\lambda^1_n}\,\left\|\OR{r}^J_n(t)\right\|_{\HL}=0;
\end{equation}
and $\widetilde{v}_n$ verifies the bound
\begin{equation}
\left\|\widetilde{v}_n\right\|_{\Snorm\left(R^d\times[-T_1t_n^{-1}\lambda^1_n,\,T_1t_n^{-1}\lambda^1_n]\right)}\leq M_1<\infty,
\end{equation}
for $J,\,n$ sufficiently large.
Consider the rescaled and translated $\widetilde{u}_n$ and $\widetilde{v}_n$ as follows:
\begin{eqnarray*}
(\widetilde{\widetilde{u}}_{0n},\widetilde{\widetilde{u}}_{1n})&:=&\left((t_n^{-1}\lambda^1_n)^{\frac{d}{2}-1}\,\widetilde{u}_n\left(t_n^{-1}\lambda^1_n\,(x+t_n^{-1}c^1_n),\,t_n^{-1}\lambda^1_n\,t\right),\right.\\
&&\quad\quad\quad\quad\quad\,\left.(t_n^{-1}\lambda^1_n)^{\frac{d}{2}}\,\partial_t\widetilde{u}_n\left(t_n^{-1}\lambda^1_n\,(x+t_n^{-1}c^1_n),\,t_n^{-1}\lambda^1_n\,t\right)\right)\\
(\widetilde{\widetilde{v}}_{0n},\widetilde{\widetilde{v}}_{1n})&:=&\left((t_n^{-1}\lambda^1_n)^{\frac{d}{2}-1}\,\widetilde{v}_n\left(t_n^{-1}\lambda^1_n\,(x+t_n^{-1}c^1_n),\,t_n^{-1}\lambda^1_n\,t\right),\right.\\
&&\quad\quad\quad\quad\quad\,\left.(t_n^{-1}\lambda^1_n)^{\frac{d}{2}}\,\partial_t\widetilde{v}_n\left(t_n^{-1}\lambda^1_n\,(x+t_n^{-1}c^1_n),\,t_n^{-1}\lambda^1_n\,t\right)\right).\\
\end{eqnarray*}
By rescaling and translation, (\ref{eq:errorvanishesinitial}) implies that for $M>1$ and $\epsilon>0$
\begin{eqnarray}
&&\left\|\left(\widetilde{\widetilde{v}}_{0n},\,\widetilde{\widetilde{v}}_{1n}\right)-\left(\widetilde{\widetilde{u}}_{0n},\,\widetilde{\widetilde{u}}_{1n}\right)\right\|_{\HL\left(B_M(0)\backslash\,\bigcup_{j\in \mathcal{J}_1\cup \mathcal{J}_2,\,j>1}B_{\epsilon}(\frac{c^j_n-c^1_n}{\lambda^1_n})\right)}\nonumber\\
&&\quad\quad\quad\quad+\,\|\widetilde{\widetilde{v}}_{0n}-\widetilde{\widetilde{u}}_{0n}\|_{L^{\dual}\left(B_M(0)\backslash\,\bigcup_{j\in \mathcal{J}_1\cup \mathcal{J}_2,\,j>1}B_{\epsilon}(\frac{c^j_n-c^1_n}{\lambda^1_n})\right)}\,\rightarrow\, 0,\label{eq:errorvanishesinitial2}
\end{eqnarray}
as $n\to\infty$.\\
Denoting 
\begin{equation}
\mathcal{J}_{\ast}:=\left\{j\in \mathcal{J}_1\cup \mathcal{J}_2:\,\,\frac{c^j_n-c^1_n}{\lambda^1_n}\,\,{\rm is\,\,bounded\,\,in\,\,}n\right\}.
\end{equation}
Passing to a subsequence, we can assume that $\frac{c^j_n-c^1_n}{\lambda^1_n}\to x_j$ as $n\to\infty$ for each $j\in \mathcal{J}_{\ast}$. For any $R>0$, $\sigma>0$, $y\in R^d$ with
\begin{equation*}
{\rm dist}\left(y,\,\{x_j,\,j\in J_{\ast}\}\right)>R+\sigma,
\end{equation*}
by Lemma \ref{lm:localinspaceapproximation} and (\ref{eq:errorvanishesinitial2}), we see that
\begin{equation}\label{eq:vanishesimportant}
\sup_{|t|<T_1\wedge R}\,\left\|\OR{\widetilde{\widetilde{u}}}_n(t)-\OR{\widetilde{\widetilde{v}}}_n(t)\right\|_{\HL(|x-y|<R-|t|)}+\left\|\widetilde{\widetilde{u}}_n-\widetilde{\widetilde{v}}_n\right\|_{\Snorm\left(\left\{(x,t):\,|x-y|<R-|t|,\,|t|<T_1\wedge R\right\}\right)}\rightarrow 0,
\end{equation}
as $n\to\infty$.\\
The vanishing condition (\ref{eq:vanishingcharacterization}), with a simple rescaling argument, implies in particular that for any $M>1$
\begin{equation*}
\lim_{n\to\infty}\,\frac{1}{t_n^{-1}\lambda^1_n}\int_{-t_n^{-1}\lambda^1_n}^{t_n^{-1}\lambda^1_n}\int_{\left|x-\frac{c^1_n}{t_n}\right|<Mt_n^{-1}\lambda^1_n}\,\left|\partial_t\widetilde{u}_n+\frac{x}{t+1}\cdot\nabla \widetilde{u}_n+\left(\frac{d}{2}-1\right)\frac{\widetilde{u}_n}{t+1}\right|^2\,dxdt=0.
\end{equation*}
Note that 
\begin{eqnarray*}
&&\frac{1}{t_n^{-1}\lambda^1_n}\int_{-t_n^{-1}\lambda^1_n}^{t_n^{-1}\lambda^1_n}\int_{\left|x-\frac{c^1_n}{t_n}\right|<Mt_n^{-1}\lambda^1_n}\,\left|\frac{\widetilde{u}_n}{t+1}\right|^2\,dxdt\\
&&\quad\lesssim \left[\sup_{t\in (-t_n^{-1}\lambda^1_n,\,t_n^{-1}\lambda^1_n)}\|\widetilde{u}_n(\cdot,t)\|^2_{L^{\dual}(R^d)}\right](t_n^{-1}\lambda^1_n)^{2}\to 0,
\end{eqnarray*}
as $n\to\infty$. Noting also that 
\begin{equation*}
\lim_{n\to\infty}\,\frac{c^1_n}{t^1_n}=\ell_1,\,\,\,{\rm and}\,\,\lambda^1_n\ll t_n,
\end{equation*}
we get that 
\begin{equation}\label{eq:vanishingcharacterization5}
\lim_{n\to\infty}\,\frac{1}{t_n^{-1}\lambda^1_n}\int_{-t_n^{-1}\lambda^1_n}^{t_n^{-1}\lambda^1_n}\int_{\left|x-\frac{c^1_n}{t_n}\right|<Mt_n^{-1}\lambda^1_n}\,\left|\partial_t\widetilde{u}_n+\ell_1\cdot\nabla \widetilde{u}_n\right|^2\,dxdt=0.
\end{equation}
Rescaling and translating (\ref{eq:vanishingcharacterization5}), we obtain that 
\begin{equation}
\lim_{n\to\infty}\int_{-1}^1\int_{|x|<M}\left|\partial_t\widetilde{\widetilde{u}}_n+\ell_1\cdot\nabla \widetilde{\widetilde{u}}_n\right|^2\,dxdt=0.
\end{equation}
By the expansion of $\widetilde{v}_n$, $\widetilde{\widetilde{v}}_n$ admits the following expansion
\begin{eqnarray*}
\widetilde{\widetilde{v}}_n&=&\OR{U}(x,t)+\OR{\widetilde{r}}^J_n\\
&&+\sum_{1\leq j\leq J,\,j\not\in \mathcal{J}_1\cup \mathcal{J}_2}\left(\left((\lambda^1_n)^{-1}\lambda^j_n\right)^{-\frac{d}{2}+1}U_j\left(\frac{x-\frac{c^j_n-c^1_n}{\lambda^1_n}}{(\lambda^1_n)^{-1}\lambda^j_n},\,\frac{t-(\lambda^1_n)^{-1}t^j_n}{(\lambda^1_n)^{-1}\lambda^j_n}\right),\right.\\
&&\quad\quad\quad\quad\quad\quad\quad\quad\quad\quad\left.\left((\lambda^1_n)^{-1}\lambda^j_n\right)^{-\frac{d}{2}}\partial_tU_j\left(\frac{x-\frac{c^j_n-c^1_n}{\lambda^1_n}}{(\lambda^1_n)^{-1}\lambda^j_n},\,\frac{t-(\lambda^1_n)^{-1}t^j_n}{(\lambda^1_n)^{-1}\lambda^j_n}\right)
\right)\\
&&\quad\quad\quad\quad\quad +\left((\lambda^1_n)^{-\frac{d}{2}+1}w^J_n\left(\frac{x+\frac{c^1_n}{\lambda^1_n}}{\lambda^1_n},\,\frac{t}{\lambda^1_n}\right),\,(\lambda^1_n)^{-\frac{d}{2}}\partial_tw^J_n\left(\frac{x+\frac{c^1_n}{\lambda^1_n}}{\lambda^1_n},\,\frac{t}{\lambda^1_n}\right)\right).
\end{eqnarray*}
By the pseudo-orthogonality condition on the parameters, we conclude that 
\begin{equation}
\OR{\widetilde{\widetilde{v}}}_n\rightharpoonup \OR{U} \,\,\,\,{\rm in}\,\,\left\{(x,t):\,|x-y|<R\wedge T_1-|t|,\,|t|<R\wedge T_1\right\}.
\end{equation}
By (\ref{eq:vanishesimportant}), we also have
\begin{equation}
\OR{\widetilde{\widetilde{u}}}_n\rightharpoonup \OR{U} \,\,\,\,{\rm in}\,\,\left\{(x,t):\,|x-y|<R\wedge T_1-|t|,\,|t|<R\wedge T_1\right\}.
\end{equation}
 Hence
\begin{equation}
\partial_tU+\ell_1\cdot\nabla U\equiv 0,
\end{equation}
in $\left\{(x,t):\,|x-y|<R\wedge T_1-|t|,\,|t|<R\wedge T_1\right\}$. \\
Moving $y$ around and shrinking $\sigma\to 0+$, we conclude that
\begin{equation}\label{eq:firstordersoliton}
\partial_tU+\ell_1\cdot\nabla U\equiv 0,
\end{equation}
in $\left\{(x,t):\,|t|<{\rm dist}\left(x,\,\{x_j,\,j\in \mathcal{J}_{\ast}\}\right)\wedge T_1\right\}$. 
(\ref{eq:firstordersoliton}) implies that $U$ satisfies
\begin{equation*}
U(x,t)=\Psi(x-\ell_1t),
\end{equation*}
in $\left\{(x,t):\,|t|<{\rm dist}\left(x,\,\{x_j,\,j\in \mathcal{J}_{\ast}\}\right)\wedge T_1\right\}$. Since $U$ is a solution to equation (\ref{eq:main}), we conclude that
\begin{equation}\label{eq:ellipticsoliton}
\ell_1\cdot\nabla\left(\ell_1\cdot\nabla \Psi\right)-\Delta\Psi=|\Psi|^{\dual-2}\Psi,
\end{equation}
in $R^d\backslash\,\{x_j,\,j\in \mathcal{J}_{\ast}\}$. By the regularity condition $\Psi\in\dot{H}^1(R^d)$, equation (\ref{eq:ellipticsoliton}) is actually satisfied across $x_j$, i.e., $x_j$ are removable singularities. Hence
\begin{equation}\label{eq:ellipticsoliton1}
\ell_1\cdot\nabla\left(\ell_1\cdot\nabla \Psi\right)-\Delta\Psi=|\Psi|^{\dual-2}\Psi,\,\,{\rm in}\,\,R^d.
\end{equation}
By Lemma 2.1 in \cite{DKMprofile}, for nontrivial $\Psi$, we must have $|\ell_1|<1$ and $\Psi\equiv Q_{\ell_1}$. Hence in this case, the profile is a travelling wave. By the energy expansion, there can be only boundedly many such profiles.\\

For profies in Case III, we have $0<\lambda^j_n\ll |t^j_n|$, then we have
\begin{eqnarray*}
&&\lim_{n\to\infty}\left\|(\lambda^j_n)^{-\frac{d}{2}+1}U_j^L\left(\frac{x-c^j_n}{\lambda^j_n},\,-\frac{t^j_n}{\lambda^j_n}\right)\right\|_{L^{\dual}(R^d)}\\
&&=\lim_{n\to\infty}\left\|U_j^L\left(x,\,-\frac{t^j_n}{\lambda^j_n}\right)\right\|_{L^{\dual}(R^d)}\\
&&=\lim_{t\to\pm\infty}\left\|U_j^L(x,t)\right\|_{L^{\dual}(R^d)}\\
&&=0,
\end{eqnarray*}
by the property of free waves. These profiles can therefore be absorbed into the error term. Recall that $v$ is the regular part of $u$, the theorem is proved.
\end{section}

\begin{section}{virial identity and exclusion of dispersive energy in the region $\{(x,t):\,|x|<\lambda t\}$ for $\lambda\in (0,1)$ along a sequence of times}
In this section, we use a virial identity to obtain a decomposition with better residue term. 
\begin{theorem}\label{th:betterresidue}
Let $u$ be as in the last section. Then there exists a time sequence $t_n\downarrow 0$, such that $\OR{u}(t_n)$ has asymptotic decomposition (\ref{eq:maindecomposition}) with properties as in Theorem \ref{th:mainpreliminary}. In addition the residue term $(\epsilon_{0n},\,\epsilon_{1n})$ satisfies the following refined vanishing condition:
\begin{equation}\label{eq:betterresidue}
\|\spartial\epsilon_{0n}\|_{L^2(B_{r_{\ast}})}+\|\nabla \epsilon_{0n}\|_{L^2(B_{\lambda t_n})}+\|\epsilon_{1n}\|_{L^2(B_{\lambda t_n})}+\left\|\epsilon_{1n}+\frac{r}{t_n}\partial_r\epsilon_{0n}\right\|_{L^2(B_{t_n})}\to 0,
\end{equation}
as $n\to\infty$, for any $\lambda\in(0,1)$.
\end{theorem}

\smallskip
\noindent
{\it Remark.} The only missing vanishing condition, in comparison with the main Theorem \ref{th:main} is that the Strichartz norm of the free radiation with initial data $(\epsilon_{0n},\,\epsilon_{1n})$ goes to zero. \\

\smallskip
\noindent
{\it Proof.} By Lemma \ref{lm:keyvanishingMorawetz}, there exists a sequence $\mu_n\downarrow 0$, such that
\begin{equation}\label{eq:vanishingcondition6}
\lim_{n\to\infty}\,\frac{1}{\mu_n}\int_{\mu_n}^{2\mu_n}\int_{|x|<4t}\,\left(\partial_tu+\frac{x}{t}\cdot\nabla u+\left(\frac{d}{2}-1\right)\frac{u}{t}\right)^2(x,t)\,dx\,dt=0,
\end{equation}
and that along two sequences of times $t_{1n}$ and $t_{2n}$ with $t_{1n}\in (\frac{4}{3}\mu_{n},\,\frac{13}{9}\mu_n)$, $t_{2n}\in (\frac{14}{9}\mu_n,\,\frac{5}{3}\mu_n)$, we have that 
\begin{equation}\label{eq:sequencevanishing1}
\sup_{0<\tau<\frac{t_{1n}}{16}}\,\frac{1}{\tau}\int_{|t-t_{1n}|<\tau}\int_{|x|<4t}\,\left(\partial_tu+\frac{x}{|x|}\cdot\nabla u+\left(\frac{d}{2}-1\right)\frac{u}{t}\right)^2\,dxdt\to 0,\,\,{\rm as}\,\,n\to\infty,
\end{equation}
and
\begin{equation}\label{eq:sequencevanishing2}
\sup_{0<\tau<\frac{t_{2n}}{16}}\,\frac{1}{\tau}\int_{|t-t_{2n}|<\tau}\int_{|x|<4t}\,\left(\partial_tu+\frac{x}{|x|}\cdot\nabla u+\left(\frac{d}{2}-1\right)\frac{u}{t}\right)^2\,dxdt\to 0,\,\,{\rm as}\,\,n\to\infty.
\end{equation}
By Theorem \ref{th:mainpreliminary}, along $t_{1n}$ and $t_{2n}$ approaching zero, we have the following decomposition
\begin{equation}\label{eq:maindecomposition5}
\OR{u}(t_{\iota n})=\OR{v}(t_{\iota n})+\sum_{j=1}^{J_{\iota}}\,\left((\lambda^j_{\iota n})^{-\frac{d}{2}+1}\, Q_{\ell_{\iota j}}\left(\frac{x-c^j_{\iota n}}{\lambda^j_{\iota n}},\,0\right),\,(\lambda^j_{\iota n})^{-\frac{d}{2}}\, \partial_tQ_{\ell_{\iota j}}\left(\frac{x-c^j_{\iota n}}{\lambda^j_{\iota n}},\,0\right)\right)+\OR{\epsilon}_{\iota n},
\end{equation}
where $\iota\in\{1,\,2\}$, and $\lambda^j_{\iota n}\ll t_{\iota n}$, $c^j_{\iota n}\in B_{\beta t_{\iota n}}$ for some $\beta\in (0,1)$ and $\ell_{\iota j}=\lim\limits_{n\to\infty}\frac{c^j_{\iota n}}{t_{\iota n}}$. $\lambda^j_{\iota n},\,c^j_{\iota n}$ satisfy the pseudo-orthogonality condition (\ref{eq:pseudo}). In the above, with slight abuse of notation, we denote $\OR{\epsilon}_{\iota n}:=(\epsilon_{\iota n},\,\partial_t\epsilon_{\iota n})\in\HL$. Moreover, $\epsilon_{\iota n}$ vanishes in the sense that for some fixed $r_0>0$
\begin{equation}\label{eq:vanishingdispersivetermLp}
\|\epsilon_{\iota n}\|_{L^{\dual}(|x|\leq r_0)}\to 0,\quad {\rm as}\,\,n\to\infty.
\end{equation}
We observe that by the decomposition (\ref{eq:maindecomposition5}), for any $\epsilon>0$ small, it holds that
\begin{eqnarray*}
\|u(\cdot,t_{\iota n})\|_{L^2(B_{t_{\iota n}})}&\lesssim& \|v(\cdot,\,t_{\iota n})\|_{L^2(|x|<t_{\iota n})}+\|\epsilon_{\iota n}\|_{L^2(|x|<t_{\iota n})}\\
                                                                 &&\quad\quad+\,\left\|\sum_{j=1}^{J_{\iota}}\,(\lambda_{\iota n})^{-\frac{d}{2}+1}\, Q_{\ell_{\iota j}}\left(\frac{x-c^j_{\iota n}}{\lambda_{\iota n}},\,0\right)\right\|_{L^2\left(\bigcup_{j=1}^{J_{\iota}}B_{\epsilon\, t_{\iota n}(c^j_{\iota n})}\right)}\\
&&\quad\quad+\,\left\|\sum_{j=1}^{J_{\iota}}\,(\lambda_{\iota n})^{-\frac{d}{2}+1}\, Q_{\ell_{\iota j}}\left(\frac{x-c^j_{\iota n}}{\lambda_{\iota n}},\,0\right)\right\|_{L^2\left(B_{t_{\iota n}}\backslash \left(\bigcup_{j=1}^{J_{\iota}}B_{\epsilon\, t_{\iota n}}(c^j_{\iota n})\right)\right)}\\
&\lesssim&\|v(\cdot,t_{\iota n})\|_{L^{\dual}(B_{t_{\iota n}})}\,t_{\iota n}+\|\epsilon_{\iota n}\|_{L^{\dual}(B_{t_{\iota n}})}\,t_{\iota n}\\
&&\quad\quad+\left\|\sum_{j=1}^{J_{\iota}}\,(\lambda_{\iota n})^{-\frac{d}{2}+1}\, Q_{\ell_{\iota j}}\left(\frac{x-c^j_{\iota n}}{\lambda_{\iota n}},\,0\right)\right\|_{L^{\dual}\left(\bigcup_{j=1}^{J_{\iota}}B_{\epsilon\, t_{\iota n}(c^j_{\iota n})}\right)}\,\epsilon\,t_{\iota n}\\
&&\quad\quad+\left\|\sum_{j=1}^{J_{\iota}}\,(\lambda_{\iota n})^{-\frac{d}{2}+1}\, Q_{\ell_{\iota j}}\left(\frac{x-c^j_{\iota n}}{\lambda_{\iota n}},\,0\right)\right\|_{L^{\dual}\left(B_{t_{\iota n}}\backslash \left(\bigcup_{j=1}^{J_{\iota}}B_{\epsilon\, t_{\iota n}}(c^j_{\iota n})\right)\right)}\,t_{\iota n}\\
&\lesssim& o(t_{\iota n})+\epsilon \,t_{\iota n}, \quad {\rm as}\,\,n\to\infty.
\end{eqnarray*}
Since this is true for all small $\epsilon>0$, we obtain that
\begin{equation}\label{eq:secondvanishing}
\lim_{n\to\infty}\,\frac{1}{t_{\iota n}}\,\|u(\cdot,t_{\iota n})\|_{L^2(B_{t_{\iota n}})}=0.
\end{equation}
Thanks to the vanishing condition (\ref{eq:secondvanishing}), we can use the following virial identity.  Multiply equation (\ref{eq:main}) with $u$, and integrate over the region $\{(x,t):\,|x|<t,\,t\in (t_{1n},\,t_{2n})\}$. We obtain
\begin{eqnarray*}
0&=&\int_{t_{1n}}^{t_{2n}}\int_{|x|<t}\,\partial_t\left(\partial_tu\,u\right)-(\partial_tu)^2-{\rm div}\,\left(\nabla u \,u\right)+|\nabla u|^2-|u|^{\dual}\,dxdt\\
&=&\frac{1}{\sqrt{2}}\,\int_{t_{1n}}^{t_{2n}}\int_{|x|=t}\,\left(-\partial_tu-\frac{x}{|x|}\cdot\nabla u\right)\,u\,d\sigma dt\\
&&\quad\quad\quad+\int_{|x|<t_{2n}}\,[\partial_tu\,u](x,t_{2n})\,dx-\int_{|x|<t_{1n}}\,[\partial_tu\,u](x,t_{1n})\,dx\\
&&\quad\quad\quad+\int_{t_{1n}}^{t_{2n}}\int_{|x|<t}\,-(\partial_tu)^2+|\nabla u|^2-|u|^{\dual}\,dxdt.
\end{eqnarray*}
Noting that $t_{2n}\sim t_{1n}$ and $t_{2n}-t_{1n}\sim t_{1n}$, with the help of the control on the energy flux, we can estimate
\begin{eqnarray*}
&&\left|\int_{t_{1n}}^{t_{2n}}\int_{|x|=t}\,\left(-\partial_tu-\frac{x}{|x|}\cdot\nabla u\right)\,u\,d\sigma dt\right|\\
&&\lesssim\,\left(\int_{t_{1n}}^{t_{2n}}\int_{|x|=t}\,\left(\partial_tu+\frac{x}{|x|}\cdot\nabla u\right)^2\,d\sigma dt\right)^{\frac{1}{2}}\,\left(\int_{t_{1n}}^{t_{2n}}\int_{|x|=t}\,u^2\,d\sigma dt\right)^{\frac{1}{2}}\\
&&\lesssim o(1)\cdot \left(\int_{t_{1n}}^{t_{2n}}\int_{|x|=t}\,|u|^{\dual}\,d\sigma dt\right)^{\frac{1}{\dual}}\,|t_{2n}-t_{1n}|\\
&&\lesssim o\left(t_{2n}-t_{1n}\right),\quad{\rm as}\,\,n\to\infty.
\end{eqnarray*}
Using crucially (\ref{eq:secondvanishing}), we can also estimate
\begin{eqnarray*}
\left|\int_{|x|<t_{\iota n}}\,\partial_tu\,u(x,t_{\iota n})\,dx\right|&\lesssim&\,\|\partial_tu\|_{L^2(B_{t_{\iota n}})}\,\left(\int_{|x|<t_{\iota n}}\,u^2(x,t_{\iota n})\,dx\right)^{\frac{1}{2}}\\
&\lesssim&o(t_{\iota n})\lesssim o(t_{2n}-t_{1n}).
\end{eqnarray*}
Summarizing the above, we conclude that
\begin{equation*}
\frac{1}{t_{2n}-t_{1n}}\int_{t_{1n}}^{t_{2n}}\int_{|x|<t}\,|\nabla u|^2-|\partial_tu|^2-|u|^{\dual}\,dxdt\to 0,\quad{\rm as}\,\,n\to\infty.
\end{equation*}
We emphasize that $t_{2n}-t_{1n}\sim \mu_n$, which is important for us.
Passing to a subsequence, we can assume that
\begin{eqnarray}
&&\left|\frac{1}{t_{2n}-t_{1n}}\int_{t_{1n}}^{t_{2n}}\int_{|x|<t}\,|\nabla u|^2-|\partial_tu|^2-|u|^{\dual}\,dxdt\right|\leq 4^{-n};\label{eq:crucialdecay2}\\
&&\frac{1}{\mu_n}\int_{\mu_n}^{2\mu_n}\int_{|x|<4t}\,\left(\partial_tu+\frac{x}{t}\cdot\nabla u+\left(\frac{d}{2}-1\right)\frac{u}{t}\right)^2(x,t)\,dx\,dt\leq 8^{-n}.
\end{eqnarray}
Denoting 
\begin{equation*}
g(t)=\int_{|x|<4t}\,\left(\partial_tu+\frac{x}{t}\cdot\nabla u+\left(\frac{d}{2}-1\right)\frac{u}{t}\right)^2(x,t)\,dx.
\end{equation*}
Then using the boundedness of Hardy-Littlewood maximal functions from $L^1$ to $L^{1,\infty}$, we get that
\begin{equation}\label{eq:smallmeasure1}
\left|\left\{t\in (t_{1n},\,t_{2n}):\,M(g\chi_{(\mu_{n},\,2\mu_{n})})(t)>2^{-n}\right\}\right|\lesssim 4^{-n}|t_{2n}-t_{1n}|.
\end{equation}
Since $\int_{R^d}\,|\nabla_{x,t}u|^2+|u|^{\dual}(x,t)\,dx\leq M$, by (\ref{eq:crucialdecay2}), we see that
\begin{equation}\label{eq:smallmeasure2}
\left|\left\{t\in (t_{1n},\,t_{2n}):\,\int_{|x|<t}\,\left[|\nabla u|^2-|\partial_tu|^2-|u|^{\dual}\right](x,t)\,dx<2^{-n}\right\}\right|\gtrsim M^{-1}2^{-n}|t_{2n}-t_{1n}|.
\end{equation}
Hence we can find a sequence $t_n\in (t_{1n},\,t_{2n})$, such that along the sequence $t_n$, we have that 
\begin{eqnarray}
&&\lim_{n\to\infty}\,M(g\chi_{(\mu_n,\,2\mu)})(t_n)=0,\label{eq:decayno1}\\
&&\limsup_{n\to\infty}\,\int_{|x|<t}\,\left[|\nabla u|^2-|\partial_tu|^2-|u|^{\dual}\right](x,t_n)\,dx\leq 0.
\label{eq:decayno2}
\end{eqnarray}
The vanishing condition (\ref{eq:decayno1}), by Theorem \ref{th:mainpreliminary},  implies that there exist $\lambda^j_{n}\ll t_n$, $c^j_n\in B_{\beta t_n}$ for some $\beta\in (0,1)$ with $\ell_j=\lim\limits_{n\to\infty}\frac{c^j_n}{t_n}$ well defined, for $1\leq j\leq J_0<\infty$, such that
\begin{equation}\label{eq:decomposition8}
\OR{u}(t_n):=\OR{v}(t_n)+\sum_{j=1}^{J_0}\left((\lambda^j_{n})^{-\frac{d}{2}+1}Q_{\ell_j}\left(\frac{x-c^j_{ n}}{\lambda^j_{n}},\,0\right),\,(\lambda^j_{n})^{-\frac{d}{2}}\partial_tQ_{\ell_j}\left(\frac{x-c^j_{ n}}{\lambda^j_{n}},\,0\right)\right)+\OR{\epsilon}_n,
\end{equation}
with $\OR{\epsilon}_n$ vanishes asymptotically in the sense that for some $r_0>0$, 
\begin{equation}\label{goodvanishing13}
\|\epsilon_n\|_{L^{\dual}(B_{r_0})}\to 0,\quad {\rm as}\,\,n\to\infty.
\end{equation}
In addition, the parameters $\lambda_n^j,\,c^j_n$ satisfy the pseudo-orthogonality condition (\ref{eq:pseudo}).\\
We now use (\ref{eq:decayno2}) to obtain refined vanishing property of $\OR{\epsilon}_n$. We remark that (\ref{eq:decayno2}) is not coercive for general solutions. However it is zero for travelling waves (see Claim \ref{claim:vanishonsoliton}) and coercive for dispersive part $\OR{\epsilon_n}$ of the solution. Hence it is particularly suited for controlling $\OR{\epsilon}_n$ along the sequence $t_n$. The idea of using additional virial type quantity such as (\ref{eq:decayno2}) to eliminate dispersive energy was firstly introduced in \cite{JiaKenig}.  We claim
\begin{claim}\label{claim:vanishonsoliton}
Let $\ell\in R^d$ with $|\ell |<1$, and let $Q_{\ell}$ be a travelling wave with velocity $\ell$. Then $Q_{\ell}(x,t)$ satisfies
\begin{equation}
\int_{R^d}\,\left[|\nabla Q_{\ell}|^2-|\partial_tQ_{\ell}|^2-|Q_{\ell}|^{\dual}\right](x,t)\,dx=0,\quad{\rm for\,\,all}\,\,t.
\end{equation}
\end{claim}
The proof will be given at the end of this section.\\

\hspace{.6cm}For the regular part $\OR{v}$, we clearly have
\begin{equation}
\lim_{t\to0+}\,\int_{|x|<t}\left[|\nabla v|^2-|\partial_tv|^2-|v|^{\dual}\right](x,t)\,dx=0.
\end{equation}
By Claim \ref{claim:vanishonsoliton}, and the pseudo-orthogonality of the profiles in the decomposition (\ref{eq:decomposition8}), we can conclude from (\ref{eq:decayno2}) that 
\begin{equation}\label{eq:decayepsilon1}
\limsup_{n\to\infty}\,\int_{|x|<t_n}\,|\nabla \epsilon_n|^2-|\partial_t\epsilon_n|^2-|\epsilon_n|^{\dual}\,dx\leq 0.
\end{equation}
By the vanishing condition $\|\epsilon_n\|_{L^{\dual}(|x|<t_n)}\to 0$, we see that
\begin{equation}
\limsup_{n\to\infty}\,\int_{|x|<t_n}|\nabla \epsilon_n|^2-|\partial_t\epsilon_n|^2\,dx=0.
\end{equation}
On the other hand, (\ref{eq:decayno1}) implies that
\begin{equation*}
\int_{|x|<t_n}\,\left(\partial_tu+\frac{x}{t}\cdot\nabla u+\left(\frac{d}{2}-1\right)\frac{u}{t}\right)^2(x,t_n)\,dx\to 0,\,\,{\rm as}\,\,n\to\infty.
\end{equation*}
By the estimates $\int_{|x|<t_n}|u|^2\,dx=o(t_n^2)$ as $n\to\infty$, which follows from similar arguments as in the proof of (\ref{eq:secondvanishing}), we also have
\begin{equation}
\int_{|x|<t_n}\,\left(\partial_tu+\frac{x}{t}\cdot\nabla u\right)^2(x,t_n)\,dx\to 0,\,\,{\rm as}\,\,n\to\infty.
\end{equation}
Using the fact that
\begin{equation}
\int_{R^d}\,\left|\partial_tQ_{\ell}+\ell\cdot\nabla Q_{\ell}\right|^2(x,t)\,dx\equiv 0,
\end{equation}
from the decomposition (\ref{eq:decomposition8}) and the orthogonality of the profiles, we conclude that
\begin{equation}\label{eq:decayepsilon2}
\lim_{n\to\infty}\,\int_{|x|<t_n}\,\left|\partial_t\epsilon_n+\frac{x}{t_n}\cdot\nabla \epsilon_n\right|^2\,dx=0.
\end{equation}
Combining (\ref{goodvanishing13}), (\ref{eq:decayepsilon1}) and (\ref{eq:decayepsilon2}), we see that
\begin{equation}\label{eq:gooddecayepsilon}
\limsup_{n\to\infty}\,\int_{|x|<t_n}\,|\nabla \epsilon_n|^2-\left|\frac{x}{t_n}\cdot\nabla \epsilon_n\right|^2\,dx\leq 0.
\end{equation}
(\ref{eq:gooddecayepsilon}) implies that for any $\lambda\in(0,1)$
\begin{eqnarray}
&&\int_{|x|<t_n}|\spartial \epsilon_n|^2\,dx\to 0,\,\,{\rm as}\,\,n\to\infty,\label{eq:epsilondecay1}\\
&&\int_{|x|<\lambda t_n}|\nabla\epsilon_n|^2\,dx\to 0,\,\,{\rm as}\,\,n\to\infty.\label{eq:epsilondecay2}
\end{eqnarray}
Combining (\ref{eq:epsilondecay2}) with (\ref{eq:decayepsilon2}), we see that in fact for any $\lambda\in(0,1)$
\begin{equation}
\int_{|x|<\lambda t_n}\,|\nabla_{x,t}\epsilon_n|^2\,dx\to 0,\,\,{\rm as}\,\,n\to\infty.
\end{equation}
The theorem is proved.\\

\smallskip
\noindent
{\it Proof of Claim \ref{claim:vanishonsoliton}.} We can assume without loss of generality that $\ell=le_1$, where $e_1=(1,0,\dots,0)\in R^d$. Then 
\begin{equation*}
Q_{\ell}=Q\left(\frac{x_1-lt}{\sqrt{1-l^2}},x_2,x_3\right).
\end{equation*}
Then direct calculations imply that
\begin{eqnarray*}
&&\int_{R^d}\left(\,|\nabla Q_{\ell}|^2-|\partial_tQ_{\ell}|^2-|Q_{\ell}|^{\dual}\right)(x,t)\,dx\\
&&=\int_{R^d}\,\left(|\nabla Q|^2-|Q|^{\dual}\right)\left(\frac{x_1-lt}{\sqrt{1-l^2}},x_2,x_3\right)\,dx\\
&&=\sqrt{1-l^2}\int_{R^d}\,\left(|\nabla Q|^2-|Q|^{\dual}\right)(x)\,dx\\
&&=0.
\end{eqnarray*}
In the last identiy we have used the fact that $Q$ is a steady state to equation (\ref{eq:main}). The claim is proved.\\

\end{section}

\begin{section}{Ruling out profiles from time infinity and proof of the main theorem}
In this section, we use the refined vanishing condition (\ref{eq:betterresidue}) to rule out any remaining nontrivial profiles, and conclude the proof of the main theorem.

\hspace{.6cm} Let us firstly prove the following result, which describes profile decomposition of a sequence of special initial data.
\begin{lemma}\label{lm:refinedprofile}
Suppose that $(u_{0n},\,u_{1n})$ is a bounded sequence in $\HL$. In addition, assume that there exists $\lambda_n\uparrow 1$, such that
\begin{equation}\label{refinedvan1}
\|u_{0n}\|_{L^{\dual}(R^d)}+\|u_{1n}\|_{L^2(R^d\backslash B_1)}+\|\spartial\, u_{0n}\|_{L^2(R^d)}\to 0,
\end{equation}
and
\begin{equation}\label{refinedvan2}
\|u_{1n}+\partial_r u_{0n}\|_{L^2(B_1)}+\|\partial_r u_{0n}\|_{L^2\left(B_{|x|<\lambda_n}\bigcup\{x:\,|x|>1\}\right)}\to 0,
\end{equation}
as $n\to\infty$. Then passing to a subsequence, $(u_{0n},\,u_{1n})$ has the following profile decomposition
\begin{eqnarray}
(u_{0n},\,u_{1n})&=&\sum_{j=1}^J\,\left(\left(\lambda^j_n\right)^{-\frac{d}{2}+1}U_j^L\left(\frac{x-c^j_n}{\lambda^j_n}, \,-\frac{t^j_n}{\lambda^j_n}\right),\,\left(\lambda^j_n\right)^{-\frac{d}{2}}\partial_tU_j^L\left(\frac{x-c^j_n}{\lambda^j_n}, \,-\frac{t^j_n}{\lambda^j_n}\right)\right)\nonumber\\
&&\nonumber\\
&&\quad\quad\quad+\left(w^J_{0n},\,w^J_{1n}\right),\label{eq:refinedprofiledecomposition}
\end{eqnarray}
where the profiles and parameters satisfy, in additional to the usual orthogonality conditions (\ref{eq:profileorthogonality1}) (\ref{eq:profileorthogonality2}), that
\begin{eqnarray}
&&\lim_{n\to\infty}\frac{t^j_n}{\lambda^j_n}\in\{\pm \infty\};\label{eq:refined1}\\
&&\lim_{n\to\infty}|t^j_n|=1;\label{eq:refined2}\\
&&\lim_{n\to\infty}c^j_n=0,\label{eq:refined3}
\end{eqnarray}
for each $j$.
\end{lemma}

\smallskip
\noindent
{\it Proof.} Passing to a subsequence, we can assume that 
$(u_{0n},\,u_{1n})$ has the following profile decomposition
\begin{eqnarray}
(u_{0n},\,u_{1n})&=&\sum_{j=1}^J\,\left(\left(\lambda^j_n\right)^{-\frac{d}{2}+1}U_j^L\left(\frac{x-c^j_n}{\lambda^j_n}, \,-\frac{t^j_n}{\lambda^j_n}\right),\,\left(\lambda^j_n\right)^{-\frac{d}{2}}\partial_tU_j^L\left(\frac{x-c^j_n}{\lambda^j_n}, \,-\frac{t^j_n}{\lambda^j_n}\right)\right)\nonumber\\
&&\nonumber\\
&&\quad\quad\quad+\left(w^J_{0n},\,w^J_{1n}\right).\label{eq:oldprofiledecomposition}
\end{eqnarray}
Let us firstly show that $\lim\limits_{n\to\infty}\frac{t^j_n}{\lambda^j_n}\in\{\pm\infty\}$. Assume for some $j$, $t^j_n\equiv 0$. By the concentration property of $(u_{0n},\,u_{1n})$, $c^j_n$ is bounded. By the assumption that $\|u_{0n}\|_{L^{\dual}(R^d)}\to 0$ as $n\to\infty$, the profile must satisfy $U^L_j(\cdot,0)\equiv 0$. Hence
\begin{equation}
\int_{R^d}|\partial_tU^L_j+\partial_rU^L_j|^2(x,0)\,dx=\int_{R^d}|\partial_tU^L_j|^2(x,0)\,dx>0.
\end{equation}
By orthogonality of profiles, we then have
\begin{equation}
\liminf_{n\to\infty}\|u_{1n}+\partial_ru_{0n}\|^2_{L^2(B_1)}\ge \left\|\left(\partial_tU^L_j+\partial_rU^L_j\right)(\cdot,0)\right\|^2_{L^2(R^d)}>0.
\end{equation}
A contradiction with (\ref{refinedvan2}). Hence there are no profiles with $t^j_n\equiv 0$.

\hspace{.6cm}For the remaining profiles, we divide into several cases.\\

{\it Case I.} Let us firstly rule out profile $\OR{U}_j^L$ with (by passing to a subsequence if necessary) $\tau=\lim\limits_{n\to\infty}|t^j_n|\neq 1$. From (\ref{eq:refined1}), we see that $\lambda^j_n\to 0$. 
A moment of reflection, using elementary geometry, shows that the set
\begin{equation}
E^j_n:=\{x:\,\lambda_n<|x|<1\}
\end{equation}
satisfies for all $M>0$ that
\begin{equation}
\left|\left\{x\in R^d:\,\left||x-c^j_n|-|t^j_n|\right|<M\lambda^j_n\right\}\cap\, E^j_n\right|=o\left(\left|t^j_n\right|^{d-1}\lambda^j_n\right),
\end{equation}
as $n\to\infty$. Thus we can apply Lemma \ref{lm:removedcone} and conclude that for any $\beta>1$ and sufficiently large $n$,
\begin{equation}
\|(u_{0n},\,u_{1n})\|_{\HL\left(\left\{x:\,\frac{|t^j_n|}{\beta}<|x-c^j_n|<\beta |t^j_n|\right\}\backslash\,E^j_n\right)}\ge \|\OR{U}^L_j\|_{\HL}>0.
\end{equation}
This is a contradiction with the property of $(u_{0n},\,u_{1n})$ that the energy is concentrated in $E^j_n$. \\

{\it Case II.} Now we rule out the profile $\OR{U}_j^L$, with (passing to a subsequence if necessary) $\lim_{n\to\infty}t^j_n\in\{\pm\infty\}$. This case is similar to {\it Case I.} We omit the routine details.\\

{\it Case III.} We have $t^j_n=o(1)$ as $n\to\infty$. We note that by the concentration property of $(u_{0n},\,u_{1n})$, $\left|c^j_n\right|\to 1$. Passing to a subsequence and by rotating the coordinate system if necessary, we can assume that $c^j_n\to e_1:=(1,0\dots,0)\in R^d$. The vanishing conditions (\ref{refinedvan1}) and (\ref{refinedvan2}) of $(u_{0n},\,u_{1n})$ imply that
\begin{equation}\label{goodv12}
\int_{B_{2t^j_n}(c^j_n)}\left|u_{1n}+\partial_1u_{0n}\right|^2+|\nabla_{x'}u_{0n}|^2\,dx\to 0,
\end{equation}
as $n\to\infty$. 
By Lemma \ref{lm:localizedorthogonality}, we have that for any $\beta>1$
\begin{eqnarray}
&&\liminf_{n\to\infty}\int_{\frac{t^j_n}{\beta}<|x-c^j_n|<\beta t^j_n}\,(u_{1n}+\partial_1u_{0n})^2+|\nabla_{x'}u_{0n}|^2\,dx\nonumber\\
&&\ge\int_{R^d}\,\left[\left(\partial_tU^L_j+\partial_1U^L_j\right)^2+|\nabla_{x'}U^L_j|^2\right](x,t)\,dx>0.\label{eq:hoholocalizedexpansion}
\end{eqnarray}
(\ref{eq:hoholocalizedexpansion}) contradicts with (\ref{goodv12}).\\

{\it Case IV.} We have $\tau=\lim\limits_{n\to\infty}|t^j_n|=1$. Using the same arguments as in {\it Case I}, we also have $\lim\limits_{n\to\infty}c^j_n=0$. These are the profiles that appear in the decomposition (\ref{eq:refinedprofiledecomposition}). \\

The lemma is proved.\\

Now we are ready to prove the main theorem.
\begin{theorem}
Let $\OR{u}\in C((0,1],\HL(R^d))$, with $u\in \Snorm(R^d\times (\epsilon,1])$ for any $\epsilon>0$, be a Type II blow up solution to equation (\ref{eq:main}) with blow up time $t=0$. Define the singular set 
\begin{equation}
\mathcal{S}:=\left\{x_{\ast}\in R^d:\,\|u\|_{\Snorm\left(B_{\epsilon}(x_{\ast})\times (0,\epsilon]\right)}=\infty,\quad{\rm for\,\,any\,\,}\epsilon>0\right\}.
\end{equation}
Then $\mathcal{S}$ is a set of finitely many points only. Assume that $0\in\mathcal{S}$ is a singular point. Then there exist integer $J_0\ge 1$, $r_0>0$, $\OR{v}\in \HL$, time sequence $t_n\downarrow 0$,  scales $\lambda_n^j$ with $0<\lambda_n^j\ll t_n$, positions $c_n^j\in R^d$ satisfying $c_n^j\in B_{\beta t_n}$ for some $\beta\in(0,1)$ with $\ell_j=\lim\limits_{n\to\infty}\frac{c_n^j}{t_n}$ well defined, and travalling waves $Q_{\ell_j}$, for $1\leq j\leq J_0$, such that inside the ball $B_{r_0}$ we have
\begin{equation}\label{eq:maindecomposition}
\OR{u}(t_n)=\OR{v}+\sum_{j=1}^{J_0}\,\left((\lambda_n^j)^{-\frac{d}{2}+1}\, Q_{\ell_j}\left(\frac{x-c_n^j}{\lambda_n^j},\,0\right),\,(\lambda_n^j)^{-\frac{d}{2}}\, \partial_tQ_{\ell_j}\left(\frac{x-c_n^j}{\lambda_n^j},\,0\right)\right)+(\epsilon_{0n},\,\epsilon_{1n}).
\end{equation}
Moreover, let $\OR{\epsilon}^L_n$ be the solution to linear wave equation with initial data $(\epsilon_{0n},\,\epsilon_{1n})\in\HL$, then the following vanishing condition holds:
\begin{equation}\label{eq:refinedintro}
\|\epsilon^L_n\|_{\Snorm(R^d\times R)}+\|( \epsilon_{0n},\,\epsilon_{1n})\|_{\HL(\{x:\,|x|<\lambda t_n,\,{\rm or}\,|x|>t_n\})}+\|\spartial\,\epsilon_{0n}\|_{L^2(R^d)}\to 0,
\end{equation}
as $n\to\infty$ for any $\lambda\in(0,1)$. In the above $\spartial$ denotes the tangential derivative.
In addition, the parameters $\lambda_n^j,\,c_n^j$ satisfy the pseudo-orthogonality condition
\begin{equation}
\frac{\lambda_n^j}{\lambda_n^{j'}}+\frac{\lambda_n^{j'}}{\lambda_n^j}+\frac{\left|c^j_n-c^{j'}_n\right|}{\lambda^j_n}\to\infty,
\end{equation}
as $n\to\infty$, for each $1\leq j\neq j'\leq J_0$.
\end{theorem}

\smallskip
\noindent
{\it Proof.} By Theorem \ref{th:betterresidue}, we only need to show that $J_0\ge 1$ and 
\begin{equation}\label{eq:lastclaim}
\lim_{n\to\infty}\|\epsilon^L_n\|_{\Snorm(R^d\times R)}=0.
\end{equation}
By Lemma \ref{lm:refinedprofile} with a rescaling, passing to a subsequence, $(\epsilon_{0n},\,\epsilon_{1n})$ has the following profile decomposition
\begin{eqnarray}
(\epsilon_{0n},\,\epsilon_{1n})&=&\sum_{j=1}^J\,\left(\left(\lambda^j_n\right)^{-\frac{d}{2}+1}U_j^L\left(\frac{x-c^j_n}{\lambda^j_n}, \,-\frac{t^j_n}{\lambda^j_n}\right),\,\left(\lambda^j_n\right)^{-\frac{d}{2}}\partial_tU_j^L\left(\frac{x-c^j_n}{\lambda^j_n}, \,-\frac{t^j_n}{\lambda^j_n}\right)\right)\nonumber\\
&&\nonumber\\
&&\quad\quad\quad+\left(w^J_{0n},\,w^J_{1n}\right),\label{eq:lastrefinedprofiledecomposition}
\end{eqnarray}
where the profiles and parameters satisfy, in additional to the usual orthogonality conditions (\ref{eq:profileorthogonality1}) (\ref{eq:profileorthogonality2}) and vanishing condition for $w^J_n$, that
\begin{eqnarray}
&&\lim_{n\to\infty}\frac{t^j_n}{\lambda^j_n}\in\{\pm \infty\};\label{eq:lastrefined1}\\
&&\lim_{n\to\infty}\frac{|t^j_n|}{t_n}=1;\label{eq:lastrefined2}\\
&&\lim_{n\to\infty}\frac{c^j_n}{t_n}=0,\label{eq:lastrefined3}
\end{eqnarray}
for each $j$. \\

We shall firstly rule out the profiles with $\lim\limits_{n\to\infty}\frac{t^j_n}{t_n}=1$. Consider the following sequence of initial data
\begin{eqnarray}\label{eq:modifiedinitialdata}
(\widetilde{u}_{0n},\,\widetilde{u}_{1n})&:=&\sum_{j=1}^J\,\left(\left(\lambda^j_n\right)^{-\frac{d}{2}+1}U_j^L\left(\frac{x-c^j_n}{\lambda^j_n}, \,-\frac{t^j_n}{\lambda^j_n}\right),\,\left(\lambda^j_n\right)^{-\frac{d}{2}}\partial_tU_j^L\left(\frac{x-c^j_n}{\lambda^j_n}, \,-\frac{t^j_n}{\lambda^j_n}\right)\right)\nonumber\\
&&\quad\quad\quad+\OR{v}(t_n)+(w^J_{0n},\,w^J_{1n}).
\end{eqnarray}
Hence we have removed any possible solitons from $\OR{u}(t_n)$, in order to be able to control the evolution of solution $\OR{\widetilde{u}}$ to equation (\ref{eq:main}) with $\OR{\widetilde{u}}(t_n)=(\widetilde{u}_{0n},\,\widetilde{u}_{1n})$. Note that
\begin{equation}\label{lastdifference}
\lim_{n\to\infty}\,\|(u_{0n},\,u_{1n})-(\widetilde{u}_{0n},\,\widetilde{u}_{1n})\|_{\HL(|x|>\frac{1+\beta}{2} t_n)}=0.
\end{equation}
Lemma \ref{lm:nonlinearprofiledecomposition} and the definition of nonlinear profiles imply that for $t\in [\frac{t_n}{2},\,t_n]$, $\OR{\widetilde{u}}$ has the following asymptotic expansion
\begin{eqnarray}
\OR{\widetilde{u}}(t)&=&\sum_{j=1}^J\,\left(\left(\lambda^j_n\right)^{-\frac{d}{2}+1}U_j^L\left(\frac{x-c^j_n}{\lambda^j_n}, \,\frac{t-t_n-t^j_n}{\lambda^j_n}\right),\,\left(\lambda^j_n\right)^{-\frac{d}{2}}\partial_tU_j^L\left(\frac{x-c^j_n}{\lambda^j_n}, \,\frac{t-t_n-t^j_n}{\lambda^j_n}\right)\right)\nonumber\\
&&\quad\quad\quad+\OR{v}(t)+\OR{w}^J_n(t)+\OR{r}^J_n(t),
\end{eqnarray}
where we recall that $\OR{w}^J_n(t)$ is the solution to the linear wave equation with $\OR{w}^J_n(t_n)=(w^J_{0n},\,w^J_{1n})$, and
\begin{equation*}
\lim_{J\to\infty}\,\limsup_{n\to\infty}\,\sup_{t\in[t_n/2,t_n]}\|\OR{r}^J_n\|_{\HL}=0.
\end{equation*}
At $t=\frac{t_n}{2}$, if $\lim\limits_{n\to\infty}\frac{t^j_n}{t_n}=1$, then by Lemma \ref{lm:concentrationoffreeradiation} and (\ref{eq:lastrefined1},\,\ref{eq:lastrefined2},\,\ref{eq:lastrefined3}), the energy of
$$\left(\left(\lambda^j_n\right)^{-\frac{d}{2}+1}U_j^L\left(\frac{x-c^j_n}{\lambda^j_n}, \,\frac{t-t_n-t^j_n}{\lambda^j_n}\right),\,\left(\lambda^j_n\right)^{-\frac{d}{2}}\partial_tU_j^L\left(\frac{x-c^j_n}{\lambda^j_n}, \,\frac{t-t_n-t^j_n}{\lambda^j_n}\right)\right)$$
is concentrated in $\{x:\,\frac{3t_n}{2\gamma}<|x|<\frac{3\gamma t_n}{2}\}$, for any $\gamma>1$. Hence by Lemma \ref{lm:localizedorthogonality}, $\OR{\widetilde{u}}(\frac{t_n}{2})$ has uniform amount of nontrivial energy accumulation in $\{x:\,\frac{3t_n}{2\gamma}<|x|<\frac{3\gamma t_n}{2}\}$. Note also by (\ref{lastdifference}) and Lemma \ref{lm:localinspaceapproximation}, we have 
\begin{equation}
\|\OR{u}(t_n)-\OR{\widetilde{u}}(t_n)\|_{\HL(\{x:\,\frac{3t_n}{2\gamma}<|x|<\frac{3\gamma t_n}{2}\})}\to 0,
\end{equation}
as $n\to\infty$ if $\gamma$ is chose close to $1$. Thus $\OR{u}(t_n)$ also has uniform amount of nontrivial energy accumulation in  $\{x:\,\frac{3t_n}{2\gamma}<|x|<\frac{3\gamma t_n}{2}\}$, as $n\to\infty$. A contradiction with the fact that outside the lightcone $|x|>t$, $\OR{u}$ is regular. \\

Next we rule out profiles with $\lim\limits_{n\to\infty}\frac{t^j_n}{t_n}=-1$, keeping in mind that the profiles with $\lim\limits_{n\to\infty}\frac{t^j_n}{t_n}=1$ have been ruled out. Introduce $\OR{\widetilde{u}}$ as in the first part of the proof. By Lemma \ref{lm:nonlinearprofiledecomposition} and the definition of nonlinear profiles, $\OR{\widetilde{u}}$ has the following asymptotic expansion for $t\in [t_n,\delta]$:
\begin{eqnarray}
\OR{\widetilde{u}}(t)&=&\sum_{j=1}^J\,\left(\left(\lambda^j_n\right)^{-\frac{d}{2}+1}U_j^L\left(\frac{x-c^j_n}{\lambda^j_n}, \,\frac{t-t_n-t^j_n}{\lambda^j_n}\right),\,\left(\lambda^j_n\right)^{-\frac{d}{2}}\partial_tU_j^L\left(\frac{x-c^j_n}{\lambda^j_n}, \,\frac{t-t_n-t^j_n}{\lambda^j_n}\right)\right)\nonumber\\
&&\quad\quad\quad+\OR{v}(t)+\OR{w}^J_n(t)+\OR{r}^J_n(t),\label{lastapproximation}
\end{eqnarray}
where $\OR{w}^J_n(t_n)=(w^J_{0n},\,w^J_{1n})$, and
\begin{equation*}
\lim_{J\to\infty}\,\limsup_{n\to\infty}\,\sup_{t\in[t_n,\delta]}\|\OR{r}^J_n\|_{\HL}=0.
\end{equation*}
Fix $\beta<\beta_1<\beta_2<1$. Define the modified sequence of initial data 
\begin{equation}
(u_{0n}',\,u_{1n}'):=\left\{\begin{array}{ll}
                                        \OR{u}(t_n)&\,\,{\rm if}\,\,|x|>\beta_2t_n,\\
   &\\
\frac{t_n^{-1}|x|-\beta_1}{\beta_2-\beta_1}\OR{u}(t_n)+\frac{\beta_2-t_n^{-1}|x|}{\beta_2-\beta_1}\OR{\widetilde{u}}(t_n)&\,\,{\rm if}\,\,\beta_1t_n<|x|<\beta_2t_n,\\
&\\
                       \OR{\widetilde{u}}(t_n) &  \,\,{\rm if}\,\,|x|<\beta_1.
                               \end{array}\right.
\end{equation}
It is easy to verify that for $r_0>0$ sufficiently small,
$$\lim_{n\to\infty}\|(u_{0n}',\,u_{1n}')-\OR{\widetilde{u}}(t_n)\|_{\HL(B_{r_0})}=0.$$
Let $\OR{u'}$ be the solution to equation (\ref{eq:main}) with $\OR{u'}(t_n)=(u_{0n}',\,u_{1n}')$. By finite speed of propagation, we see that
\begin{equation}\label{eq:identical}
u(x,t)=u'(x,t),\,\,{\rm for}\,\,|x|>\beta_2\,t_n+t-t_n,\,t\in [t_n,\delta].
\end{equation}
By Lemma \ref{lm:longtimeperturbation},  and by choosing $\delta$ smaller if necessary, we see that 
\begin{equation}\label{thesame}
\sup_{t\in[t_n,\,\delta]}\,\|\OR{u'}(t)-\OR{\widetilde{u}}(t)\|_{\HL(B_{r_0-\delta})}\to 0,\,\,{\rm as}\,\,n\to\infty.
\end{equation}
By (\ref{lastapproximation}) and Lemma \ref{lm:localizedorthogonality}, at time $t=\delta$, $$\liminf_{n\to\infty}\left\|\OR{\widetilde{u}}(\delta)\right\|_{\HL(\delta+\beta_2t_n-t_n<|x|<\delta+t_n)}\ge\|\OR{U}^L_j\|_{\HL(R^d)}.$$
Hence by (\ref{thesame}) and (\ref{eq:identical}), $\OR{u}(\delta)$ has uniform nontrivial energy concentration $\delta+\beta_2t_n-t_n<|x|<\delta+t_n$, as $n\to\infty$. A contradicition with the fact that $\OR{u}(\delta)\in\HL$. 

\hspace{.6cm}In summary, there are no nontrivial profiles in the decomposition (\ref{eq:lastrefinedprofiledecomposition}). Thus, $\epsilon^L$ verifies (\ref{eq:lastclaim}). Then $J_0\ge 1$ follows since $0\in \mathcal{S}$. The theorem is proved.\\

\medskip

\center{{\bf Acknowledgement}\\
We thank Carlos Kenig for the valuable discussions during the preparation of this work.}

\medskip
\end{section}

\end{document}